\documentclass[a4paper,11pt]{amsart}
\usepackage{amsmath,inputenc,euscript,amssymb,geometry}
\usepackage{graphicx}
\usepackage{amssymb}
\usepackage{latexsym}
\usepackage{amssymb,amsbsy,amsmath,amsfonts,amssymb,amscd}
\usepackage{hyperref}

\newtheorem{lemma}{Lemma}[section]
\newtheorem{theorem}[lemma]{Theorem}
\newtheorem{corollary}[lemma]{Corollary}
\newtheorem{remark}[lemma]{Remark}

\newcommand{\R}{\mathbb{R}}
\newcommand{\D}{\displaystyle}

\title{On the stability of a singular vortex dynamics}
\thanks{E-mail addresses : Valeria.Banica@univ-evry.fr, mtpvegol@lg.ehu.es}

\begin{document}

\maketitle
\begin{center}{{\bf{V. Banica$^1$, L. Vega$^{2,*}$}}  \vspace{3mm}\\\tiny{$^1$D\'epartement de Math\'ematiques, Universit\'e d'Evry, France\\
$^2$Departamento de Matemáticas, Universidad del Pais Vasco, Spain }}
\end{center}
\begin{abstract}In this paper we address the question of the singular vortex dynamics exhibited in \cite{SJL}, which generates a corner in finite time. The purpose is  to prove that under some appropriate small regular perturbation  the corner still remains. Our approach uses the Hasimoto transform and deals with the long range scattering properties of a Gross-Pitaevski equation with time-variable coefficients. \end{abstract}

\tableofcontents

\section{Introduction}
In this paper we study the stability properties of selfsimilar solutions of the geometric flow
\begin{equation}\label{binormal}
\chi_t=\chi_x\land\chi_{xx}.
\end{equation}
Here $\chi=\chi(t,x)\in\mathbb R^3$, $x$ denotes the arclength parameter and $t$ the time variable. The above equation was proposed by DaRios in 1906 \cite{DaR} and re-derived by Arms and Hama in 1965 \cite{ArHa} as an approximation of the dynamics of a vortex filament under Euler equations. In this model $\chi(t,x)$ represents the support of the singular vectorial measure that describes the vorticity. The velocity field is then obtained from the Biot-Savart integral and is singular at the points of the filament. Equation \eqref{binormal} follows from a Taylor expansion around a given point. The first term is  discarded by symmetry, and then, after doing a re-normalization in time to avoid a logarithmic singularity, the second term gives \eqref {binormal}.  Therefore just local effects are considered and for this reason this model is usually known as the Localized Induction Approximation (LIA). We refer the reader to \cite{Ba} and \cite{Sa} for an analysis and discussion about the limitations of this model and to \cite{Ri} for a survey about Da Rios' work. 

Starting with the work by Schwartz in \cite{Sch} LIA has been also used as an approximation of the quantum vortex motion in superfluid Helium. In particular in the  recent work by T. Lipniacki \cite{Li1}, \cite{Li2}, a detailed analysis of the selfsimilar solutions of \eqref{binormal} is also made. A rather complete list of references about the use of LIA in this setting can be found in these two papers.\par

Let us recall now that for a general curve in $\mathbb{R}^3$, parametrized by arclength, its tangent vector $T$, its normal vector $n$ and its binormal $b$ satisfy the Frenet system
\begin{equation}\label{Frenet}
\left(\begin{array}{c}
T\\n\\b
\end{array}\right)_x=
\left(\begin{array}{ccc}
0 & c & 0 \\ -c & 0 & \tau \\ 0 &  -\tau & 0 
\end{array}\right)
\left(\begin{array}{c}
T\\n\\b
\end{array}\right),
\end{equation}
where $c$ is the curvature of the curve and $\tau$ its torsion.
Then equation \eqref{binormal} can be rewritten as 
\begin{equation}\label{binormal1}
\chi_t=cb.
\end{equation}
This explains why the term binormal flow is sometimes used as a substitute to LIA.

Another relevant connection of \eqref{binormal} is obtained by computing the equation satisfied by the tangent vector 
$$T=\chi_x.$$
An immediate calculation gives that $T$ has to solve 
\begin{equation}\label{map}
T_t=T\land T_{xx}.
\end{equation}
Notice that as a consequence the arclength parametrization is preserved and therefore $T$ gives a flow onto the unit sphere $\mathbb{S}^2$. Equation \eqref{map} can be rewritten as
\begin{equation}\label{magn}
T_t=JD_xT_x,
\end{equation}
with $J$ denoting the complex structure of the sphere and $D_x$ the covariant derivative. With this formulation we identify \eqref{magn} as the Schr\"odinger map onto the sphere. This equation can be also seen as a simplification of Landau-Lifchitz equation for ferromagnetism (see \cite{LaLi}).\par

$ $

As we have already said our main interest is to study the selfsimilar solutions of the binormal flow \eqref{binormal}. Let us recall the known results about selfsimilar solutions. Although there is a one parameter family of possible scalings that leave invariant the set of solutions, there is only one which preserves the property that $T(t,x)$ is in $\mathbb{S}^2$. Namely, for $\lambda>0$, if $\chi(t,x)$ solves \eqref{binormal}, so does $$\frac1\lambda\,\chi(\lambda^2 t,\lambda x).$$ 
Let us look for solutions of the type
$$\chi(t,x)=\sqrt{t}\,G\left(\frac {x}{\sqrt{t}}\right).$$
After differentiation we get that $G$ has to be a solution of the ODE
$$\frac 12G-\frac x2G'=G'\land G''.$$
Computing another derivative and with some abuse of notation  we get that $T(x)$ has to solve
\begin{equation}\label{eqTself}
-\frac x2 T'=T\land T''=(cb)'.
\end{equation}
Using the Frenet equations \eqref{Frenet} it follows that
$$-\frac x2 cn=c'b-c\tau n.$$
As a conclusion we obtain a one parameter family of curves (see \cite{La1}, \cite{La2}, \cite{Bu}) characterized by 
\begin{equation}\label{ctauselfsim}
c(x)=a\qquad,\qquad\tau(x)=\frac{x}{2}.
\end{equation}
Let us notice that \eqref{eqTself} implies
\begin{equation}\label{idTself}
\left(T+\frac{2c}{x}b\right)'=-\frac{2c}{x^2}b.
\end{equation}We define $(T_a,n_a,b_a)(x)$ to be the unique solution of the Frenet system with curvature and torsion as in \eqref{ctauselfsim}, and initial data $(T_a,n_a,b_a)(0)=I_3$. 
By using the fact that the binormal vector is unitary, 
an immediate consequence of \eqref{idTself} is that that for any $a\in\mathbb R^+$ there exists a pair of unit vectors $A^\pm_a$ such that
\begin{equation}\label{T_a}
\lim_{x\rightarrow\pm\infty}T_a(x)=A^\pm_a.
\end{equation}
In \cite{SJL} it is proved, among other things, the following result.

{\bf{Theorem}} (Gutierrez-Rivas-Vega){\it{ Let $a$ be a positive number, and let $G_a$ be defined by $G'_a=T_a$ with $G_a(0)=2a(0,0,1)$. Then\\

(i)$\qquad \qquad\qquad\left|\sqrt{t}G_a\left(\frac{x}{\sqrt{t}}\right)-A^+_a x \mathbb{I}_{[0,\infty)}(x)-A^-_a x \mathbb{I}_{(-\infty,0]}(x)\right|\leq {2a}{\sqrt{t}}.$\\

(ii) For any test function $\psi(x)$ such that 
$$\int_{-\infty}^\infty|\psi(x)|\frac{dx}{1+|x|}<\infty,$$
we have}}
\begin{equation}\label{ii}
\lim_{t\rightarrow 0}\int_{-\infty}^\infty \left(T_a\left(\frac{x}{\sqrt{t}}\right)-A^+_a \mathbb{I}_{[0,\infty)}(x)-A^-_a \mathbb{I}_{(-\infty,0]}(x)\right)\psi(x)dx=0.
\end{equation}
\noindent
{\it{
(iii) The relation between $a$ and $A^\pm_a$ is
$$\sin\frac\theta 2=e^{-\frac{a^2}{2}},$$
where $\theta$ is the angle between the vectors $A^+_a$ and $-A^-_a$.}}
\\

Notice that \eqref{binormal} is invariant under rotations. 
As a conclusion, there exists a solution $\chi_a$ of \eqref{binormal} with the initial condition
$$\chi_a(0,x)=A^+ x \mathbb{I}_{[0,\infty)}(x)+A^- x \mathbb{I}_{(-\infty,0]}(x),$$
for any pair of unit vectors $A^\pm$, different and non-opposite.
This is deduced first by determining the number $a$ such that (iii) holds for $A^\pm$, then taking $\chi_a(x,t)=\sqrt tG_a(x/\sqrt t)$ with $G_a$ given in the Theorem,  and finally by applying to $\chi_a$ the rotation that sends $A^\pm$ into $A^\pm_a$.

Also notice that  \eqref{binormal}  is a time reversible flow because if $\chi(t,x)$ is a solution, so is $\chi(-t,-x)$. Therefore if we look at \eqref{binormal}  backwards in time  with the initial condition at  $t=1$ given by
$$\chi_a(1,x)=G_a(x),$$
we get an example of a solution which is regular at $t=1$, in fact real analytic, and that develops a singularity in the shape of a corner at time zero. 

The main result of this paper is given in Theorem \ref{stab} where we prove that under a smallness assumptions on $a$, there exist regular solutions $\chi$ of \eqref{binormal} for $t>0$, perturbations of $\chi_a$, that still have a corner at $t=0$.

\begin{remark}
Equation \eqref{magn} suggests many possible generalizations by considering other targets besides the sphere. 
Let us then introduce the notation
$$u\land_\pm v=\left(\begin{array}{ccr}
1&0&0\\
0&1&0\\
0&0&\pm1
\end{array}\right)u\land v.$$
Therefore instead of \eqref{magn} we write
\begin{equation}\label{mapsign}
T_t=T\land_\pm T_{xx},
\end{equation}
with $T$ a map from $\R^2$ onto the sphere $\mathbb{S}^2$ or the hyperbolic plane $\mathbb{H}^2$ depending on which sign is considered in \eqref{mapsign}. The positive sign stands for $\mathbb{S}^2$, and the negative one for $\mathbb{H}^2$. Analogously, since  $T=\chi_x$, we can obtain the equation 
\begin{equation}\label{binormalsign}
\chi_t=\chi_x\land_\pm\chi_{xx}.
\end{equation}
Similar calculations and results as in \cite{SJL} were done by de la Hoz \cite{Pa} for selfsimilar solutions of \eqref{binormalsign} in the hyperbolic setting. The extra-difficulty is that there is a-priori no control on the size of the euclidean length of the generalized binormal vectors. \end{remark}

In order to write our results we need to recall another remarkable connection of \eqref{binormalsign} made by Hasimoto in \cite{Ha}. It is as follows. Assume that $\chi$ is a regular solution of \eqref{binormalsign} with a strictly  positive curvature at all points. He defines the "filament function" as
\begin{equation}\label{Hasimoto}
u(t,x)=c(t,x)\exp{\left\{i\int\limits_0^x \tau(t,x')dx'\right\}}.
\end{equation}
Then $u$ solves the nonlinear Schr\"odinger equation  
\begin{equation}\label{Hasimoto-eq}
iu_t+u_{xx}\pm\frac{1}{2}\left(|u|^2-A(t)\right)u=0,
\end{equation}
with
\begin{equation}\label{Hasimoto-coeff}
A(t)=\left(\pm2\frac{c_{xx}-c\,\tau^2}{c}+c^2\right)(t,0).
\end{equation}
Let us notice that the identity \eqref{Hasimoto-coeff} provides us of some extra information on  $c$ and $\tau$ at $x=0$. This will be an important ingredient in the proof of Theorem \ref{stab}. 

As we see the focusing  sign ($+$) is related to the sphere, while the defocusing  sign (-) is connected to hyperbolic space. The real coefficient $A(t)$ can be easily eliminated by an integrating factor so that \eqref{Hasimoto-eq} can be reduced to the well known cubic NLS. This equation is completely integrable and among the infinitely many conserved quantities it has, we want to recall that
\begin{equation}\label{mass}\int|u(t,x)|^2dx\end{equation}
is preserved. This quantity is related to the kinetic energy of the filament (see \cite{Ri2}). \\

The particular selfsimilar solution $\chi_a$ has as a filament function
$$u_a(t,x)=a\frac{e^{i\frac{x^2}{4t}}}{\sqrt{t}}.$$
Therefore neither $u_a$ nor any of its derivatives are in $L^2$. As a consequence none of the other conserved quantities are finite for $u_a$, included \eqref{mass}. However, 
 we will see below \eqref{energy} that there is  a natural energy asociated to $u_a$.
 
 Notice  that $u_a$ is a solution of equation \eqref{Hasimoto-eq} with 
$$A(t)=\frac{a^2}{t},$$
and 
$$u_a(0,x)=a\delta_0.$$
Therefore, in order to study perturbations of the particular solution $u_a$ we have to study \eqref{Hasimoto-eq} within a functional setting which includes functions of infinite energy. This was started in \cite{VaVe} and then extended in \cite{Gr} and in \cite{Ch} to the case of periodic boundary conditions. None of these works consider initial data as singular as the delta-function which is our interest here. There is an obstruction to do that, as it was observed in \cite{KPV}. Using the so-called Galilean invariance the authors proved that the solution of \eqref{Hasimoto-eq} with $A(t)=0$ and $u(0,x)=a\delta_0$ either does not exist or is not unique. The reason is that the natural candidate for such a solution is 
$$\frac{a}{\sqrt{t}}e^{\pm ia^2\log t+i\frac{x^2}{4t}},$$
which has no limit as $t$ goes to zero. As noticed, in our case the Hasimoto transform leads to \eqref{Hasimoto-eq}-\eqref{Hasimoto-coeff} with $A(t)=\frac{a^2}{t}$ and therefore we have the solution $u_a$ even at $t=0$. \\

The study of the stability of $u_a$ was started in \cite{VL} where a weak stability result is obtained. We proceed as follows. First notice that after a rescaling, \eqref{Hasimoto-eq} with $A(t)= a^2/t$ can be rewritten as
\begin{equation}\label{vortex}
iu_t+u_{xx}\pm\left(|u|^2-\frac{a^2}{t}\right)u=0.
\end{equation}
Consider $u$ a solution of \eqref{vortex} for any $x\in\R$ and $t>0$. Using the so called  pseudo-conformal transformation we define a new unknown $v$ as
\begin{equation}\label{calT}
u(t,x)=\mathcal {T}v(t,x)=\D\frac{e^{i\frac{x^2}{4t}}}{\sqrt{t}}\overline{v}\left(\frac 1t,\frac xt\right).
\end{equation}
Then $\,v\,$ solves 
\begin{equation}\label{GP}
iv_t+v_{xx}\pm\D\frac 1t\left(|v|^2-a^2\right)v=0,\\
\end{equation}
and $v_a=a$ is a particular solution. A natural quantity associated to \eqref{GP} is the normalized energy
\begin{equation}\label{energy}E(t)=\frac{1}{2}\int |v_x(t)|^2\,dx\mp\frac{1}{4t}\int(|v(t)|^2-a^2)^2\,dx.
\end{equation}
An immediate calculation gives that
$$\partial_{t}E(t)\mp\frac{1}{4t^2}\int(|v|^2-a^2)^2\,dx=0.$$
In \cite{VL} we use this energy law to prove that \eqref{GP} is globally (respectively locally) well posed for $t>t_0>0$ if $E(v(t_0,x))<\infty$ in the defocusing (respectively focusing) settings. The global existence follows by proving the control
$$\|v(t)-a\|_{L^2}< C \sqrt{t}.$$
Let us notice here that similar tools have been also used by Tsutsumi and Yajima in \cite{TsYa} to prove scattering in $L^2$ for NLS with $H^1\cap\Sigma$ data.
\\

Our first theorem can be seen as an extension of the results in \cite{VL}. We construct modified wave operators for $v-a$ in both focusing and defocusing cases, under some smallness assumptions. Since we are working around a non-integrable particular solution, the source term of equation \eqref{GP} is the linear one, with a coefficient with decay $\frac{1}{t}$, that is exactly the frame for long range effects for cubic 1-d NLS (\cite{Ozawa},\cite{Ca},\cite{HaNa}). Here the situation is different since the $L^\infty$-norm of the functions we are working with is not decaying as $t$ goes to infinity, being just bounded. A link can be made with the  scattering for the Gross-Pitaevskii equation around the constant solution (\cite{GNT}), but still the situation is pretty different. 

Given $u_+$ we define 
$$v_1(t,x)=a+e^{\pm i a^2\log t}e^{it\partial_x^2}u_+(x).$$

\begin{theorem}\label{longrange}
Let $t_0>0$. 
There exists a constant $a_0>0$ such that for all $a<a_0$, and for all $u_+$ small in $L^1\cap L^2$ with respect to $a_0$ and to $t_0$, the equation (\ref{GP}) has a unique solution 
$$v-v_1\in \mathcal{C}([t_0,\infty),L^2(\mathbb{R}))\cap L^4([t_0,\infty),L^\infty(\mathbb{R})),$$ 
verifying, as $t$ goes to infinity,
\begin{equation}\label{rate}
\|v(t)-v_1(t)\|_{L^2}+\|v-v_1\|_{L^4((t,\infty),L^\infty)}=\mathcal{O}(t^{-\frac 14}).
\end{equation}
\end{theorem}

Let us notice that the family of solutions we have found is such that 
$$\|v(t)-a\|_{L^2}=\mathcal{O}(1),$$
as $t$ goes to infinity, while as we said before, for a general solution of \eqref{GP}, with sign -, we got in \cite{VL} only a control in $\mathcal{O}(\sqrt{t})$.

Once $v$ is obtained, we recover $u$ by the pseudo-conformal transformation \eqref{calT}. 
If we define
$$u_1(t,x)=a\frac{e^{i\frac{x^2}{4t}}}{\sqrt{t}}+\frac{e^{\pm ia^2\log t}}{\sqrt{4\pi i}}\hat {\overline{u}}_+\left(-\frac x2\right),$$
we get the following corollary from Theorem \ref{longrange}.

\begin{corollary}\label{thvortex}
Let $\tilde{t_0}>0$. 
There exists a constant $a_0>0$ such that for all $a<a_0$, for all $u_+$ small in $L^1\cap L^2$ with respect to $a_0$ and to $T_0$, $u_+$ in $L^2(x^4dx)$, the equation (\ref{vortex}) has a unique solution 
$$u-u_1\in \mathcal{C}((0,\tilde{t_0}],L^2(\mathbb{R}))\cap L^4((0,\tilde{t_0}],L^\infty(\mathbb{R})),$$ 
verifying, as $t$ goes to zero,
\begin{equation}\label{rateu}
\|u(t)-u_1(t)\|_{L^2}+\|u-u_1\|_{L^4((0,t),L^\infty)}=\mathcal{O}(t^{\frac 14}).
\end{equation}
In particular,
\begin{equation}\label{modulusu}
\left\| \left|u(t,x)-a\frac{e^{i\frac{x^2}{4t}}}{\sqrt{t}}\right|^2-\left|\hat{\overline{u}}_+\left(-\frac x2\right)\right|^2\right\|_{L^1}=\mathcal{O}(t^{\frac 14}),
\end{equation}
\begin{equation}\label{modulusu2}
\left\| u(t,x)-a\frac{e^{i\frac{x^2}{4t}}}{\sqrt{t}}\right\|_{L^2}=\mathcal{O}(1),
\end{equation}
but there is no limit in $L^2$ for $u(t,x)-a\frac{e^{i\frac{x^2}{4t}}}{\sqrt{t}}$ as $t$ goes to zero.
\end{corollary}
 We shall prove below that although $u$ does not have a trace at $t=0$ we will be able to construct a family of curves associated to $u$ by the Hasimoto transform that do have a limit at $t=0$.

The proof of Theorem \ref{longrange} goes as follows. We write
$$v-a=w.$$
If $v$ solves \eqref{GP} then $w$ solves
$$iw_t+w_{xx}\pm\D\frac 1t\left(|a+w|^2-a^2\right)(a+w)=0.$$
The source term includes therefore two linear terms, namely 
$$\frac{a^2}{t}w,\qquad and \qquad\frac{a^2}{t}\overline{w}.$$
As a first guess we can treat $w$ as a perturbation from a free evolution at $t=\infty$. Therefore let us assume that for $t$ large,
$$w(t)\approx e^{it\partial_x^2}u_+,$$
where $e^{it\partial_x^2}u_+$ denotes the solution of the free Schr\"odinger equation with $u_+$ as initial condition.
Then the two linear terms lead to the Duhamel integrals
$$a^2\int_t^\infty e^{i(t-\tau)\partial_x^2}e^{i\tau\partial_x^2}u_+\frac{d\tau}{\tau},$$
and 
\begin{equation}\label{Duh2}
a^2\int_t^\infty e^{i(t-\tau)\partial_x^2}e^{-i\tau\partial_x^2}\overline{u}_+\frac{d\tau}{\tau}.
\end{equation}
Clearly there is no cancellation in the first integral which therefore diverges. As a conclusion, the initial ansatz has to be modified to 
$$(v-a)e^{\mp ia^2\log t}=w.$$
Doing this the second integral \eqref{Duh2} still remains (in fact a harmless variation of it). But in this case plenty of cancellations can be expected. We exploit them by the so-called Strichartz estimates \cite{St} (see beginning of \S \ref{Lpsection}). Notice that we are in one dimension case and these estimates were proved by Fefferman and Stein in \cite {Fe}. For our later purposes the rate of decay of $v-v_1$ is crucial. The power $t^{-\frac 14}$ is proved using the mixed norm spaces $L^4 L^\infty-L^{4/3}L^1$ introduced by Ginibre and Velo \cite{GiVe} and cannot be improved if we use standard Strichartz estimates.

A natural question  is how to construct the curves $\chi(t,x)$ from the solutions obtained in Corollary \ref{thvortex}. There could be a problem if we want to use the Frenet frame, because we do not know if $|u|\neq 0$, and so the torsion cannot be well-defined by \eqref{Hasimoto}. This can be overcomed by using another type of frames \cite{Ko}. In \cite{NaShVeZe} it is proved how to construct $\chi(t,x)$ for $t>0$ and therefore to solve \eqref{binormalsign} with  regularity assumptions similar to those given by Corollary \ref{thvortex}.  The necessary modifications are straightforward. However the existence of a trace at $t=0$ is very unclear. Moreover the
$-\frac 14$ rate of decay doesn't seem enough in order to prove that the formation of a corner is preserved. \\

The main content of our second result is the improvement of the rate of decay of $v-v_1$ by strengthening the conditions on $u_+$. As we said, the first test to be checked is to obtain a better rate of convergence for the oscillatory integral \eqref{Duh2}. We are going to proceed in a different way. Recall that if $\hat u_+$ denotes the Fourier transform of $u$ we obtain the identity
$$e^{it\partial_x^2}u_+=\int e^{-it\xi^2+ix\xi}\hat{u}_+(\xi)d\xi.$$
Plugging this in the integral \eqref{Duh2} we get after changing the order of integration
$$\int_t^\infty e^{i(t-\tau)\partial_x^2}e^{-i\tau\partial_x^2}\overline{u}_+\frac{d\tau}{\tau}=\int e^{ix\xi}\hat{\overline{u}}_+(\xi)\int_t^\infty e^{i(t-2\tau)\xi^2}\frac{d\tau}{\tau}d\xi,$$
where the last integral has to be understood as an oscillatory one
\begin{equation}\label{intosc}\lim_{R\rightarrow \infty}\int_t^R e^{i\tau\xi^2}\frac{d\tau}{\tau}=\frac{e^{it\xi^2}}{t\xi^2} + remainder.
\end{equation}
This suggests to consider data $u_+\in\dot{H}^{-2}$ which we define as
\begin{equation}\label{condu}
\frac{\hat{u}_+}{|\xi|^2}\in L^2.
\end{equation}
Similar conditions were assumed by Bourgain and Wang in \cite{BoWa}.\par
For $s, p\in\mathbb{N}^*$ ,  $W^{s,p}$ is defined as
$$W^{s,p}=\{f \,|\, \nabla^kf\in L^p,\, \forall\, 0\leq k\leq s\},$$
and  $H^s=W^{s,2}.$ 
We have the following theorem.

\begin{theorem}\label{longrangeH}
Let $t_0>0$,  $s\in\mathbb{N}^*$. 
There exists a constant $a_0>0$ such that for all $a<a_0$, for all $u_+$ small in $\dot{H}^{-2}\cap H^{s}\cap W^{s,1}$ with respect to $a_0$ and to $t_0$, the equation (\ref{GP}) has a unique solution
$$v-v_1\in \mathcal{C}([t_0,\infty),H^s(\mathbb{R})),$$
verifying, as $t$ goes to infinity, and for all integer $0<k\leq s$,
\begin{equation}\label{rateH}\|(v-v_1)(t)\|_{L^2}=\mathcal{O}(t^{-\frac 12})\,\,\,,\,\,\,\|\nabla^k(v-v_1)(t)\|_{L^2}=\mathcal{O}(t^{-1}).\end{equation}
\end{theorem}
Notice that from \eqref{intosc}-\eqref{condu} we expect a $1/t$ decay coming from the linear term. However, in Theorem \ref{longrangeH} we  obtain just
$1/\sqrt{t}$ as rate of the $L^2$ convergence. The problem comes now from the quadratic terms. They are the following ones 
$$\frac{2a}{t}|w|^2,\qquad and\qquad \frac{a}{t}w^2.$$
Again the first one gives less cancellations than the second one. However the derivative $\partial_x |w|^2$ behaves better and a rate of decay $1/t$ is also proved in this case (see Lemma \ref{quadraticosc}). Similar ideas have been used in \cite{GNT}.\\

Theorem \ref{longrangeH} is enough for our purposes. First notice that by taking $u_+$ small and regular enough, we get from  Theorem \ref{longrangeH} a  solution $v$ regular and not vanishing. Hence we can define (\S\ref{ctau}) a regular curvature and torsion by taking respectively the modulus and the derivative of the phase of $u=\mathcal{T}(v)$. 
Then, we can use the Frenet frame to construct (\S\ref{constrsection}) a family of curves $\chi(t,x)$ that is a solution of the binormal flow for $t>0$. To be able to use the Frenet frame is particularly useful for us, because for example equation \eqref{idTself} has a natural generalization, see \eqref{idT}, which plays an important role in proving that the constructed family of curves is close to the selfsimilar one $\chi_a$. Also, for proving this final fact, the strong rates of decay of Theorem \ref{longrangeH} are crucial.

Finally let us say that in the construction of the family of curves we shall deal just with solutions of \eqref{binormal} and not of \eqref{binormalsign}. The obstruction for doing it in the second case is the same as the one mentioned before: there is no a-priori control on the size of the euclidean length of the generalized binormal vector if we work in $\mathbb{H}^2$. This does not happen in the sphere setting where we are able to prove the existence of the trace of $\chi$ at $t=0$ from the uniform bound of the curvature obtained in \S\ref{ctau} from Theorem \ref{longrangeH},
$$|c(t,x)|< \frac{C}{\sqrt{t}}.$$
Using this bound in \eqref{binormal1} together with the fact that  $b$ is unitary, we get the integrability of $\chi_t$ at $t=0$, and therefore  the existence of a curve $\chi_0(x)=\chi(0,x)$ follows immediately.

Our final result is the following one. 
 
\begin{theorem}\label{stab}
We fix $\epsilon>0$, $\tilde{t_0}>0$ and a positive number $a$ such that $a<a_0$, where $a_0$ is the constant in Theorem \ref{longrangeH}. Let $u_+$ small enough in $\dot{H}^{-2}\cap H^3\cap W^{3,1}$ with $x^2u_+$ small enough in $H^1$ in terms of $\epsilon$, $a$ and $\tilde{t_0}$, and let $v$ be the corresponding solution obtained in Theorem \ref{longrangeH}. By using the Hasimoto transform, we construct from $v$ a family of curves $\chi(t,x)$ which solves for $\tilde{t_0}>t>0$
$$\chi_t=c\, b,$$
and such that there exists a unique $\chi_0$ such that 
$$|\chi(t,x)-\chi_0(x)|<Ca\sqrt t$$
uniformly on $x\in(-\infty,\infty)$.\\

Moreover $\chi_0$ is Lipschitz and for $x>0$
$$|\chi_0(x)-\chi_0(0)-A^+_ax|<\epsilon \,x,$$
and
$$|\chi_0(0)-\chi_0(-x)-A^-_ax|<\epsilon \,x,$$
with  $A^+_a$ and $A^-_a$ the vectors given in \eqref{T_a} and that satisfy 
$$\sin\frac\theta 2=e^{-\frac{a^2}{2}},$$
where $\theta$ is the angle between $A^+_a$ and $-A^-_a$.
\end{theorem}

In the proof we show that the tangent vector of the binormal flow $\chi$ we construct is close to the one of $\chi_a$. We then prove that $\chi$ is close to $\chi_a$ even at time $t=0$. We recall here that
$$\chi_a(0,x)=A^+_a x \mathbb{I}_{[0,\infty)}(x)+A^-_a x \mathbb{I}_{(-\infty,0]}(x).$$

As a conclusion of the statement of Theorem \ref{stab} we get that  $\chi_0(x)$ lies in the $\epsilon$-cone around $\chi_a(0,x)$, and therefore a corner is still formed for $\chi(t,x)$ at $t=0$ and at $x=0$. The angle of this corner can be made as close as desired to the one in between $A^+_a$ and $-A^-_a$  by taking $\epsilon$ small enough.
The family of perturbations of $\chi_a$ that we obtain is determined by the wave operator constructed in Theorem \ref{longrangeH}. A better description of the allowed perturbations would be obtained if the asymptotic completeness of this wave operator were proved. This will be done in a forthcoming paper. \\

The paper is organized as follows.
In the next section we give the results about the wave operator. We start by writing the  long-range profile that implies a modification of the free evolution, and then we find the corresponding integral equation associated to this profile. In the next subsections, \S\ref{Lpsection} and \S\ref{Hsection}, we solve this integral equation first in the mixed norm spaces (Theorem \ref{longrangeH}), and then in the Sobolev spaces (Theorem \ref{longrange}). Subsection \S\ref{Corsection} is devoted to the proof of Corollary \ref{thvortex}. 
 
 Section \S\ref{flowsection} contains the construction of the family of curves $\chi$ associated to the solution  $v$ obtained in  Theorem \ref{longrangeH}, and that solve \eqref{binormal}. First, in \S\ref{estvsection} we obtain estimates on $v$ from Theorem \ref{longrangeH}. In \S\ref{ctau}, after defining the curvature and the torsion from $v$, we compute their leading terms as $t$ goes to zero. With this curvature and torsion, we construct in \S\ref{constrsection} a binormal flow up to $t=0$, as stated in the first part of  Theorem \ref{stab}. 
 
 Section \S\ref{formationsection} is devoted to the proof of the fact that the constructed flow $\chi$ is close to $\chi_a$. In the three first subsections we show that the tangent vector of $\chi$ is close to the one of $\chi_a$, and in \S\ref{cornersection} we conclude the second part of Theorem \ref{stab}.

In the last section we derive some extra-information on $\chi(t,x)$. Finally in the appendix we sketch how to construct the tangent, normal and binormal vectors of a solution of \eqref{binormal} from a solution of \eqref{Hasimoto-eq}.\\

Acknowledgements: First author was partially supported by the ANR project ``\'Etude qualitative des EDP". Second author was partially supported by the grant MTM 2007-62186 of MEC (Spain) and FEDER . Part of this work was done while the second author was visiting the University of Cergy-Pontoise.

\section{Modified wave operators}
\bigskip First we give the fixed point argument that we use to obtain the wave operator for our problem. Subsection \S\ref{Lpsection} contains the proof of Theorem \ref{longrange} in mixed norm spaces, and subsection \S\ref{Hsection} deals with the proof of Theorem \ref{longrangeH} in the Sobolev space framework. In subsection \S\ref{Corsection} we prove Corollary \ref{thvortex}.

As usual for nonlinear Schr\"odinger equations, if we want a solution of equation \eqref{GP} to behave as $t$ goes to infinity  like a particular function $v_1$, it is enough to find a fixed point for the operator
$$Av(t)=v_1(t)+i\int_t^{\infty} e^{i(t-\tau)\partial_x^2}\left(\mp\frac{(|v|^2-a^2)v}{\tau}-(i\partial_\tau+\partial_x^2)v_1(\tau)\right)d\tau,$$
in a space defined around $v_1$. 
We take as an ansatz for our problem
$$v_1=a+e^{it\partial_x^2}\omega,$$
with
$$\omega(t,\cdot)=u_+(\cdot)e^{ i\gamma\log t},$$
and $\gamma$ to be chosen later. It follows that
$$(i\partial_t+\partial_x^2)v_1=e^{it\partial_x^2}i\partial_t\omega=-\frac{\gamma}{t}e^{it\partial_x^2}\omega.$$
So for $v$ in some suitable space defined around $v_1$, we shall have to estimate
$$Av-v_1=i\int_t^{\infty} e^{i(t-\tau)\partial_x^2}\left(\mp\left(\frac{(|v|^2-a^2)v}{\tau}-\frac{(|v_1|^2-a^2)v_1}{\tau}\right)+\left(\mp\frac{(|v_1|^2-a^2)v_1}{\tau}+\frac{\gamma}{\tau}e^{i\tau\partial_x^2}\omega\right)\right)d\tau.$$
The first term of the right hand side will be easier to treat than the last one. We compute
$$(|v_1|^2-a^2)v_1=(v_1)^2\overline{v_1}-a^2v_1=\left(a^2+2ae^{it\partial_x^2}\omega+(e^{it\partial_x^2}\omega)^2\right)\left(a+\overline{e^{it\partial_x^2}\omega}\right)-a^2(a+e^{it\partial_x^2}\omega)$$
$$=a^2\overline{e^{it\partial_x^2}\omega}+a^2e^{it\partial_x^2}\omega+2a|e^{it\partial_x^2}\omega|^2+a(e^{it\partial_x^2}\omega)^2+|e^{it\partial_x^2}\omega|^2e^{it\partial_x^2}\omega.$$
Here we make the choice $\gamma=\pm a^2$, to get rid of one of the linear terms. By doing this, the only linear term left is out of resonance and the integral will converge. 

In conclusion, we are choosing 
$$v_1=a+e^{it\partial_x^2}\omega,$$
with
$$\omega(t,\cdot)=u_+(\cdot)e^{ \pm ia^2\log t},$$
and we shall do a fixed point argument in spaces defined around $v_1$, for the operator
\begin{equation}\label{fixedpoint}
Av=v_1\mp i\int_t^{\infty} e^{i(t-\tau)\partial_x^2}\left(\frac{(|v|^2-a^2)v}{\tau}-\frac{(|v_1|^2-a^2)v_1}{\tau}\right)d\tau\end{equation}
$$\mp i\int_t^{\infty} e^{i(t-\tau)\partial_x^2}\frac{a^2\overline{e^{i\tau\partial_x^2}\omega}+2a |e^{i\tau\partial_x^2}\omega|^2+a(e^{i\tau\partial_x^2}\omega)^2+|e^{i\tau\partial_x^2}\omega|^2e^{i\tau\partial_x^2}\omega}{\tau}d\tau.$$

\bigskip

Let us finally recall the 1-D Strichartz estimates that will be used throughout this section (see \cite{St}, \cite{GiVe}). We have
\begin{equation}\label{Str1} \left\| e^{it\partial_x^2}f \right\|_{L^{p_1}({\R};L^{q_1})}\leq C\,\|f\|_{L^2},
\end{equation}
and the inhomogeneous version ($1/r':=1-1/r$)
\begin{equation}\label{Str2} \left\| \int_{I\cap\{s\leq t\}} e^{i(t-s)\partial_x^2}F(s)ds\right\|_{L^{p_1}(I,L^{q_1})}\leq C\,\|F\|_{L^{p'_2}(I,L^{q'_2})},
\end{equation}
for any admissible couples $(p_i,q_i)$, that is
$$\frac {2}{p_i} +\frac {1}{q_i}=\frac 12,\qquad p\geq 2.$$
The admissible couples we shall use here are $(\infty,2)$ and $(4,\infty)$. Also, let us recall the dispersion inequality
\begin{equation}\label{disp}|e^{it\partial_x^2}f|\leq \frac{C}{\sqrt{t}}\|f\|_{L^1}.
\end{equation}
In particular,
\begin{equation}\label{dispv1}|v_1(t)|\leq a+C\,\frac{\|u_+\|_{L^1}}{\sqrt{t}}.\end{equation}

\subsection{Modified wave operators in mixed norm spaces: proof of Theorem \ref{longrange}}\label{Lpsection}

Let $t_0>0$. We shall perform the fixed point argument for the operator \eqref{fixedpoint} in the closed ball
$$X_R=\left\{v\,|\,\|v\|_X=\underset{t\in[t_0,\infty[}{\sup}\,t^\nu\|v(t)-v_1(t)\|_{L^2}+t^\nu\|v-v_1\|_{L^4((t,\infty)L^\infty)}\leq R\right\},$$
with $0<\nu$ and $R$ to be precised later.
Let us notice that in view of \eqref{dispv1}, a function $v\in X_R$ satisfies 
\begin{equation}\label{modulus}
\|v(t)\|_{L^\infty}\leq \|v_1(t)\|_{L^\infty}+\|v(t)-v_1(t)\|_{L^\infty}\leq a+C\,\frac{\|u_+\|_{L^1}}{\sqrt{t}}+\|v(t)-v_1(t)\|_{L^\infty}.
\end{equation}

We want, for a $v\in X_R$, to estimate in $X_R$
$$Av-v_1=\mp i\int_t^{\infty} e^{i(t-\tau)\partial_x^2}\left(\frac{(|v|^2-a^2)v}{\tau}-\frac{(|v_1|^2-a^2)v_1}{\tau}\right)d\tau$$
$$\mp i\int_t^{\infty} e^{i(t-\tau)\partial_x^2}\frac{a^2\overline{e^{i\tau\partial_x^2}\omega}+2a |e^{i\tau\partial_x^2}\omega|^2+a(e^{i\tau\partial_x^2}\omega)^2+|e^{i\tau\partial_x^2}\omega|^2e^{i\tau\partial_x^2}\omega}{\tau}d\tau=I+J_3+J_2+J_1.$$
We denote here $I$ to be the first term in the right-hand side, and $J_k$ to be the parts of the second term involving $k$-powers of $e^{i\tau\partial_x^2}\omega$.

For $I$ we shall use the inhomogeneous Strichartz estimates \eqref{Str2},
$$\|I\|_X=\left\|\int_t^{\infty} e^{i(t-\tau)\partial_x^2}\left(\frac{(|v|^2-a^2)v}{\tau}-\frac{(|v_1|^2-a)^2v_1}{\tau}\right)d\tau\right\|_X $$
$$\leq C\,\underset{t_0\leq t}{\sup}\,t^\nu\int_t^{\infty}\||v|^2v-|v_1|^2v_1-a^2(v-v_1)\|_{L^2}\frac{d\tau}{\tau}$$
$$\leq C\,\underset{t_0\leq t}{\sup}\,t^\nu\int_t^{\infty}(a^2+\|v_1\|^2_{L^\infty}+\|v\|^2_{L^\infty})\|v-v_1\|_{L^2}\frac{d\tau}{\tau}.$$
Since $v$ is in $X_R$, 
$$\|I\|_X\leq C\,\|v\|_X\,\underset{t_0\leq t}{\sup}\,t^\nu\int_t^{\infty}(a^2+\|v_1\|^2_{L^\infty}+\|v\|^2_{L^\infty}) \frac{d\tau}{\tau^{1+\nu}},$$
and by using \eqref{dispv1} and \eqref{modulus},
$$\|I\|_X\leq C\,\|v\|_X\left(a^2+\frac{\|u_+\|_{L^1}^2}{t_0}+\underset{t_0\leq t}{\sup}\,t^\nu\int_t^{\infty} \|(v-v_1)(\tau)\|^2_{L^\infty} \frac{d\tau}{\tau^{1+\nu}}\right).$$
In the last integral we apply Cauchy-Schwarz inequality to recover the $L^4L^\infty$ norm, and finally,
$$\|I\|_X\leq C\,\|v\|_X\left(a^2+\frac{\|u_+\|_{L^1}^2}{t_0}+\frac{\|v\|_X^2}{\sqrt{t_0}}\right).$$

The contribution of the cubic power of $e^{-it\partial_x^2}\omega$ is easy to estimate. By using the inhomogeneous Strichartz estimates \eqref{Str2} and the dispersion inequality \eqref{disp} we get 
$$\|J_3\|_X=\left\|\int_t^{\infty} e^{i(t-\tau)\partial_x^2}\frac{|e^{i\tau\partial_x^2}\omega|^2e^{i\tau\partial_x^2}\omega}{\tau}d\tau\right\|_X\leq C\,\underset{t_0\leq t}{\sup}\,t^\nu\int_t^\infty\||e^{i\tau\partial_x^2}\omega|^2e^{i\tau\partial_x^2}\omega|\|_{L^2}\frac{d\tau}{\tau}$$
$$\leq C\,\underset{t_0\leq t}{\sup}\,t^\nu\int_t^{\infty}\|e^{i\tau\partial_x^2}\omega\|^2_{L^\infty}\|e^{i\tau\partial_x^2}\omega\|_{L^2}\frac{d\tau}{\tau}\leq C\,\underset{t_0\leq t}{\sup}\,t^\nu\int_t^{\infty}\|u_+\|^2_{L^1}\|u_+\|_{L^2}\frac{d\tau}{\tau^2}\leq C(u_+)\,\underset{t_0\leq t}{\sup}\,\frac{t^\nu}{t}.$$

The quadratic terms can be handled in the same way, and we obtain
$$\|J_2\|_X=\left\|\int_t^{\infty} e^{i(t-\tau)\partial_x^2}\frac{2a|e^{i\tau\partial_x^2}\omega|^2}{\tau}d\tau\right\|_X\leq Ca\,\|u_+\|_{L^1}\|u_+\|_{L^2}\,\underset{t_0\leq t}{\sup}\,\frac{t^\nu}{t^\frac{1}{2}}.$$

So at the end we need to estimate only the linear term
$$J_1=\int_t^{\infty} e^{i(t-\tau)\partial_x^2}\frac{a^2\overline{e^{i\tau\partial_x^2}\omega}}{\tau}d\tau=\int_t^{\infty} e^{i(t-2\tau)\partial_x^2}\,\left(\overline{u_+}\,\frac{a^2}{\tau^{1\pm ia^2}}\right)d\tau.$$
First we estimate its $L^2$ norm in space. We use the conservation of the mass for the linear evolution \eqref{Str1}, 
$$\left\|J_1(t)\right\|_{L^2}= \left\|\int_t^{\infty} e^{-i2\tau\partial_x^2}\left(\overline{u_+}\frac{a^2}{\tau^{1\pm ia^2}}\right)d\tau\right\|_{L^2},$$
and the inhomogeneous Strichartz estimates \eqref{Str2},
$$\|J_1(t)\|_{L^2}\leq Ca^2\,\left\|\frac{\overline{u_+}}{\tau}\right\|_{L^{{p'}}((t,\infty),L^{{q'}})}=Ca^2\,\|u_+\|_{L^{{q'}}}\left\|\frac{1}{\tau}\right\|_{L^{{p'}}(t,\infty)}=Ca^2\,\|u_+\|_{L^{{q'}}}\frac{1}{t^\frac{1}{p}}.$$
Therefore
$$\underset{t_0\leq t}{\sup}\,t^\nu\left\|J_1(t)\right\|_{L^2}\leq Ca^2\,\|u_+\|_{L^{{q'}}}\,\underset{t_0\leq t}{\sup}\,\frac{t^\nu}{t^\frac{1}{p}}.$$
We need then $u_+\in L^{{q'}}$, $\nu<\frac{1}{p}$. 
From the admissibility relation, the best choice is $p=4$. 
Moreover, since $u_+\in L^1\cap L^2$, then by interpolation we have also $u_+\in L^{{q'}}$. \\
For estimating the $L^\infty$ norm in space of $J_1$, we use the dispersion inequality \eqref{disp}, 
$$\left\|J_1(t)\right\|_{L^\infty}\leq Ca^2\,\int_t^{\infty}\frac{\|u_+\|_{L^1}}{\tau(2\tau-t)^{\frac{1}{2}}}d\tau\leq Ca^2\,\|u_+\|_{L^1}\frac{1}{t^\frac{1}{2}}.$$
Then 
$$\underset{t_0\leq t}{\sup}\,t^\nu\left\|J_1\right\|_{L^4((t,\infty)L^\infty)}\leq Ca^2\,\|u_+\|_{L^1}\underset{t_0\leq t}{\sup}\,\frac{t^\nu}{t^\frac{1}{4}}.$$

In conclusion, for $u_+\in L^1\cap L^2$, we have obtained
$$\|A(v)\|_X\leq \,\|v\|_X\left(a^2+\frac{\|u_+\|_{L^1}^2}{t_0}+\frac{\|v\|_X^2}{\sqrt{t_0}}\right)+C(u_+)\,\underset{t_0\leq t}{\sup}\,\frac{t^\nu}{t}+C(u_+)\,a\,\underset{t_0\leq t}{\sup}\,\frac{t^\nu}{t^\frac{1}{2}}+
C(u_+)\,a^2\,\underset{t_0\leq t}{\sup}\,\frac{t^\nu}{t^\frac{1}{4}}.$$
Let $\nu=\frac 14$. Then there exists a constant $R$ small with respect to $t_0$, and a positive small constant $a_0$, such that for all $a\leq a_0$ and all $u_+$ small enough in $L^1\cap L^2$ with respect to $t_0,a_0$ and $R$, we can apply the fixed point theorem in $X_R$. We get then a unique solution $v$ of equation (\ref{GP}) such that
$$v-v_1\in \mathcal{C}([t_0,\infty),L^2(\mathbb{R}))\cap L^4([t_0,\infty),L^\infty(\mathbb{R})),$$ 
with the rate of decay (\ref{rate}), so the proof of Theorem \ref{longrange} is complete.

\subsection{Proof of Corollary \ref{thvortex}}\label{Corsection}
Let $\tilde{t_0}>0$. We denote $t_0=\frac {1}{\tilde{t_0}}$ and we consider $v$ to be the corresponding solution of Theorem \ref{longrange}, satisfying the decay \eqref{rate} as $t$ goes to infinity,
$$\|v-v_1\|_{L^\infty((t,\infty),L^2)\cap L^4((t,\infty),L^\infty)}=\mathcal{O}(t^{-\frac 14}).$$
Then $u$, the pseudo-conformal transform of $v$, will satisfy equation \eqref{vortex}. We want to show the first assertion \eqref{rateu} of Corollary \ref{thvortex}, namely the decay as $t$ goes to zero,
$$\|u-u_1\|_{L^\infty((0,t),L^2)\cap L^4((0,t),L^\infty)}=\mathcal{O}(t^{\frac 14}).$$
The mixed normed spaces we are using spaces are invariant under the pseudo-conformal transformation $\mathcal{T}$, and since
$$u=\mathcal{T}(v)\qquad,\qquad u_1=\mathcal{T}\left(a+e^{\pm i a^2\log t}\frac{e^{i\frac{x^2}{4t}}}{\sqrt{4\pi it}}\hat u_+\left(\frac{x}{2t}\right)\right),$$
we notice that \eqref{rateu} is equivalent to have, as $t$ goes to infinity,
$$\left\|v(t,x)-a-e^{\pm i a^2\log t}\frac{e^{i\frac{x^2}{4t}}}{\sqrt{4\pi it}}\hat u_+\left(\frac{x}{2t}\right)\right\|_{L^\infty((t,\infty),L^2)\cap  L^4((t,\infty),L^\infty)}=\mathcal{O}(t^{-\frac 14}).$$
In view of \eqref{rate}, this is equivalent to have this decay for the difference
$$\left\|v_1(t,x)-a-e^{\pm i a^2\log t}\frac{e^{i\frac{x^2}{4t}}}{\sqrt{4\pi it}}\hat u_+\left(\frac{x}{2t}\right)\right\|_{L^\infty((t,\infty),L^2)\cap  L^4((t,\infty),L^\infty)}=\mathcal{O}(t^{-\frac 14}).$$
From the definition of $v_1$, it is enough to prove
$$\left\|e^{it\partial_x^2}u_+-\frac{e^{i\frac{x^2}{4t}}}{\sqrt{4\pi it}}\hat u_+\left(\frac{x}{2t}\right)\right\|_{L^\infty((t,\infty),L^2)\cap  L^4((t,\infty),L^\infty)}=\mathcal{O}(t^{-\frac 14}),$$
which is one of the properties of the linear Schr\"odinger evolutions.
On the one hand, in $L^4((t,\infty),L^\infty)$ both terms decay like $t^{-\frac{1}{4}}$ as $t$ goes to infinity. 
On the other hand, the expression of the free Schr\"odinger solution gives
$$\left\|e^{it\partial_x^2}u_+-\frac{e^{i\frac{x^2}{4t}}}{\sqrt{4\pi it}}\hat u_+\left(\frac{x}{2t}\right)\right\|_{L^2}=
\frac{1}{\sqrt{4\pi t}}\left\|\int e^{-i\frac{xy}{2t}}e^{i\frac{y^2}{4t}}u_+(y)dy-\hat u_+\left(\frac{x}{2t}\right)\right\|_{L^2}$$
$$=\frac{1}{\sqrt{4\pi t}}\left\|\widehat{e^{i\frac{\cdot^2}{4t}}u_+(\cdot)}\left(\frac {x}{2t}\right)-\hat{u_+}\left(\frac {x}{2t}\right)\right\|_{L^2}= c\,
\left\| \left(e^{i\frac{y^2}{4t}}-1\right)u_+(y)\right\|_{L^2}.$$
If $u_+$ is in $L^2\cap L^2(y^4dy)$, that is if $\hat u_+\in H^2$, then this difference is $\mathcal{O}(t^{-1})$.
In conclusion, the first part \eqref{rateu} of Corollary \ref{thvortex} is proved. 
Relation \eqref{modulusu} is obtained from \eqref{rateu} by using the general formula
$$\||f|^2-|g|^2\|_{L^1}\leq (\|f\|_{L^2}+\|g\|_{L^2})\|f-g\|_{L^2},$$
and then \eqref{modulusu2} follows by the triangle inequality.

\subsection{Modified wave operators in Sobolev spaces: proof of Theorem \ref{longrangeH}}\label{Hsection}
Let $t_0>0$, $s\in\mathbb{N}^*$. In this subsection we shall perform the fixed point argument for the operator \eqref{fixedpoint} in the closed ball
$$Y_R=\left\{v\,|\,\|v\|_Y=\underset{t\in[t_0,\infty[}{\sup}\,|t|^\nu\|(v-v_1)(t)\|_{L^2}+\underset{1\leq k\leq s}{\Sigma}\,\,\underset{t\in[t_0,\infty[}{\sup}\,|t|^\mu\|\nabla^k(v-v_1)(t)\|_{L^2}\leq R\right\},$$
for strictly positive $\nu$, $\mu$ and $R$, to be precised later. 
Let us notice that  in one dimension  $|f|^2\leq \|f\|_{L^2}\|f'\|_{L^2}$. Then, for $v\in Y_R$,
\begin{equation}\label{estmodulusdiff2}
|(v-v_1)(t)|\leq\frac{C\,R}{t^\frac{\mu+\nu}{2}},
\end{equation}
and
\begin{equation}\label{estmodulusdiff}
|\nabla^k(v-v_1)(t)|\leq C\,\|\nabla^{k+1}(v-v_1)(t))\|_{L^2}^\frac 12\|\nabla^{k}(v-v_1)(t))\|_{L^2}^\frac 12\leq \frac{C\,R}{t^\mu},
\end{equation}
for all $0< k<s$. Moreover, by using the dispersion inequality \eqref{disp}, for all $0\leq k\leq s$,
$$\left|\nabla^ke^{\pm ia^2\log t}e^{it\partial_x^2}u_+\right|\leq C\frac{\|\nabla^k u_+\|_{L^1}}{\sqrt{t}}.$$
It follows that for $v\in Y_R$,
\begin{equation}\label{estmodulusv}
|v(t)|\leq |v_1(t)|+|v(t)-v_1(t)|\leq a+C\frac{\|u_+\|_{L^1}}{\sqrt{t}}+\frac{C\,R}{t^\frac{\mu+\nu}{2}},
\end{equation}
and that for all $0<k\leq s$,
\begin{equation}\label{estmodulusder}
|\nabla^kv(t)|\leq C|\nabla^kv_1(t)|+|\nabla^k(v-v_1)(t)|\leq C\frac{\|\nabla^ku_+\|_{L^1}}{\sqrt{t}}+\frac{C\,R}{t^\mu}.
\end{equation}

The proof follows as in the subsection \S\ref{Lpsection}, by estimating the terms $I$ and $J_k$ in $Y$. 

By using the conservation of the $L^2$ norm of the free equation \eqref{Str1}, 
$$\underset{t_0\leq t}{\sup}\,t^\nu\|I(t)\|_{L^2}=\underset{t_0\leq t}{\sup}\,t^\nu\left\|\int_t^{\infty} e^{i(t-\tau)\partial_x^2}\left(\frac{(|v|^2-a^2)v}{\tau}-\frac{(|v_1|^2-a)^2v_1}{\tau}\right)d\tau\right\|_{L^2}$$
$$\leq C\,\underset{t_0\leq t}{\sup}\,t^\nu\int_t^{\infty}\||v|^2v-|v_1|^2v_1-a^2(v-v_1)\|_{L^2}\frac{d\tau}{\tau}$$
$$\leq C\,\underset{t_0\leq t}{\sup}\,t^\nu\int_t^{\infty}(a^2+\|v_1\|^2_{L^\infty}+\|v\|^2_{L^\infty})\|v-v_1\|_{L^2}\frac{d\tau}{\tau}.$$
Since $v$ is in $Y_R$, 
$$\underset{t_0\leq t}{\sup}\,t^\nu\|I(t)\|_{L^2}\leq C\,\|v\|_Y\,\underset{t_0\leq t}{\sup}\,t^\nu\int_t^{\infty}(a^2+\|v_1\|^2_{L^\infty}+\|v\|^2_{L^\infty}) \frac{d\tau}{\tau^{1+\nu}}.$$
By using the bound \eqref{dispv1} on $v_1$, and \eqref{estmodulusv} on $v$, we get
$$\underset{t_0\leq t}{\sup}\,t^\nu\|I(t)\|_{L^2}\leq C\,\|v\|_Y\left(a^2+\frac{\|u_+\|_{L^1}^2}{t_0}+\frac{(C\,R)^2}{t_0^{\mu+\nu}}\right).$$
For the $L^2$ norm of the first derivative, we have
$$\underset{t_0\leq t}{\sup}\,t^\mu\|\nabla I(t)\|_{L^2}\leq C\,\underset{t_0\leq t}{\sup}\,t^\mu\int_t^{\infty}\|\nabla(|v|^2v-|v_1|^2v_1)\|_{L^2}+a^2\|\nabla(v-v_1)\|_{L^2}\frac{d\tau}{\tau}$$
$$\leq Ca^2\|v\|_Y+C\,\underset{t_0\leq t}{\sup}\,t^\mu\int_t^{\infty}(\|v\|^2_{L^\infty}\|\nabla(v-v_1)\|_{L^2}+(\|v\|_{L^\infty}+\|v_1\|_{L^\infty})\|\nabla v_1\|_{L^\infty}\|v-v_1\|_{L^2})\frac{d\tau}{\tau}.$$
By using again the fact that $v\in Y_R$ and the bounds \eqref{dispv1}, \eqref{estmodulusv}, we obtain
$$\underset{t_0\leq t}{\sup}\,t^\mu\|\nabla I(t)\|_{L^2}\leq C\,\|v\|_Y\left(a^2+\frac{\|u_+\|_{L^1}^2}{t_0}+\frac{(C\,R)^2}{t_0^{\mu+\nu}}\right)\left(1+\underset{t_0\leq t}{\sup}\,\frac{\,t^\mu}{t^{\frac{1}{2}+\nu}}\right).$$
The higher order derivatives can be estimated similarly, and we get
$$\|I\|_{Y}\leq C\,\|v\|_Y\left(a^2+\frac{\|u_+\|_{L^1}^2}{t_0}+\frac{(C\,R)^2}{t_0^{\mu+\nu}}\right)\left(1+\underset{t_0\leq t}{\sup}\,\frac{\,t^\mu}{t^{\frac{1}{2}+\nu}}\right).$$

We estimate now $J_3$ by using the invariance of the $H^s$ norm for the free evolution, and by using the fact that $H^s(\mathbb R)$ is an algebra, 
$$\|J_3(t)\|_{H^s}=\left\|\int_t^{\infty} e^{i(t-\tau)\partial_x^2}\frac{|e^{i\tau\partial_x^2}\omega|^2e^{i\tau\partial_x^2}\omega}{\tau}d\tau\right\|_{H^s}\leq \,\underset{t_0\leq t}{\sup}\,t^\mu\int_t^\infty\||e^{i\tau\partial_x^2}\omega|^2e^{i\tau\partial_x^2}\omega|\|_{H^s}\frac{d\tau}{\tau}$$
$$\leq C\,\int_t^{\infty}\|e^{i\tau\partial_x^2}\omega\|^2_{L^\infty}\|e^{i\tau\partial_x^2}\omega\|_{H^s}\frac{d\tau}{\tau}
\leq \frac{C}{t}\,\|u_+\|^2_{L^1}\|u_+\|_{H^s}.$$ 

We first consider the $L^2$ norm of $J_2$, that can be estimated as done in the previous subsection \S\ref{Lpsection},
$$\underset{t_0\leq t}{\sup}\,t^\nu\|J_2(t)\|_{L^2}\leq C\,a\,\underset{t_0\leq t}{\sup}\,t^\nu\int_t^{\infty}\||e^{i\tau\partial_x^2}\omega|^2|\|_{L^2}\frac{d\tau}{\tau}$$
$$\leq C\,a\,\underset{t_0\leq t}{\sup}\,t^\nu\int_t^{\infty}\|e^{i\tau\partial_x^2}\omega\|_{L^\infty}\|e^{i\tau\partial_x^2}\omega\|_{L^2}\frac{d\tau}{\tau}\leq C\,a\,\|u_+\|_{L^1}\|u_+\|_{L^2}\,\underset{t_0\leq t}{\sup}\,t^\nu\int_t^{\infty}\frac{d\tau}{\tau^\frac 32}$$
$$=C\,a\,\|u_+\|_{L^1}\|u_+\|_{L^2}\,\underset{t_0\leq t}{\sup}\,\frac{t^\nu}{t^\frac 12}.$$
Of course, the derivatives can also be estimated in this way. Nevertheless, for our final purpose of studying the binormal flow, we shall need more decay on the derivatives. More precisely, we have the following lemma concerning $J_2(t)$. 

\begin{lemma}\label{quadraticosc}
If $u_+\in \dot{H}^{-1}\cap \dot{H}^{s-1}$, then
\begin{equation}\label{quadraticHs}
\Sigma_{0<k\leq s}\left\|\nabla^k J_2(t)\right\|_2\lesssim\frac{a}{t}\left(\|u_+\|_{\dot{H}^{-1}}^2+\|u_+\|^2_{\dot{H}^{s-1}}\right).
\end{equation}
\end{lemma}

\begin{proof}
We have
$$\nabla^kJ_2=\nabla^ka\int_t^{\infty} e^{i(t-\tau)\partial_x^2}\frac{\left(e^{i\tau\partial_x^2}\omega\right)^2}{\tau}d\tau+2\nabla^ka\int_t^{\infty} e^{i(t-\tau)\partial_x^2}\frac{\left|e^{i\tau\partial_x^2}\omega\right|^2}{\tau}d\tau.$$

By using the Fourier transform in space, the first quadratic term is
$$\nabla^k\int_t^{\infty} e^{i(t-\tau)\partial_x^2}\frac{\left(e^{i\tau\partial_x^2}\omega\right)^2}{\tau}d\tau=\int_t^{\infty} \int\frac{e^{-i(t-\tau)\xi^2}}{\tau^{1\mp 2ia^2}}\,|\xi|^k\,e^{ix\xi}\widehat{e^{i\tau\partial_x^2}u_+}* \widehat{e^{i\tau\partial_x^2}u_+}\,d\xi\, d\tau$$
$$=\int_t^{\infty} \int\int\frac{e^{-it\xi^2-i\tau(\eta^2+(\eta-\xi)^2-\xi^2)}}{\tau^{1\mp 2ia^2}}\,|\xi|^k\,e^{ix\xi}\hat{u_+}(\eta) \hat{u_+}(\eta-\xi)\,d\eta\,d\xi \,d\tau.$$
We perform an integration by parts in time and get
$$\nabla^k\int_t^{\infty} e^{i(t-\tau)\partial_x^2}\frac{\left(e^{i\tau\partial_x^2}\omega\right)^2}{\tau}d\tau=- \int\int\frac{e^{-it(\xi^2+2\eta(\eta-\xi))}}{t^{1\mp 2ia^2}}\,|\xi|^k\,e^{ix\xi}\,\frac{\hat{u_+}(\eta) \hat{u_+}(\eta-\xi)}{2i\eta(\eta-\xi)}\,d\eta\,d\xi $$
$$-(1\mp 2ia^2)\int_t^{\infty}\int\int\frac{e^{-it\xi^2-2i\tau\eta(\eta-\xi)}}{\tau^{2\mp 2ia^2}}\,|\xi|^k\,e^{ix\xi}\,\frac{\hat{u_+}(\eta) \hat{u_+}(\eta-\xi)}{2i\eta(\eta-\xi)}\,d\eta\,d\xi\,d\tau.$$
Therefore, as $t$ goes to infinity,
$$\left\|\nabla^ka\int_t^{\infty} e^{i(t-\tau)\partial_x^2}\frac{\left(e^{i\tau\partial_x^2}\omega\right)^2}{\tau}d\tau\right\|_2\lesssim \frac{a}{t}\left\||\xi|^k\left(\frac{\hat{u_+}}{\cdot}*\frac{\hat u_+}{\cdot}\right)(\xi)\right\|_{L^2}.$$
For $k\geq 1$, $H^k$ is an algebra, so we get 
$$
\Sigma_{0<k\leq s}\left\|\nabla^ka\int_t^{\infty} e^{i(t-\tau)\partial_x^2}\frac{\left(e^{i\tau\partial_x^2}\omega\right)^2}{\tau}d\tau\right\|_2\lesssim \frac{a}{t}\left\|\mathcal{F}\left(\frac{\hat{u_+}}{\cdot}\right)\right\|_{H^k}^2\lesssim \frac{a}{t}\left(\|u_+\|_{\dot{H}^{-1}}^2+\|u_+\|_{\dot{H}^{k-1}}^2\right).
$$

The second quadratic term can be also estimated similarly and we get
$$\nabla^k\int_t^{\infty} e^{i(t-\tau)\partial_x^2}\frac{\left|e^{i\tau\partial_x^2}\omega\right|^2}{\tau}d\tau=\int_t^{\infty}\int|\xi|^k\frac{e^{-i(t-\tau)\xi^2}}{\tau}\,e^{ix\xi}\widehat{e^{i\tau\partial_x^2}u_+}* \widehat{e^{-i\tau\partial_x^2}\overline{u}_+}\,d\xi\, d\tau$$
$$=\int_t^{\infty} \int\int\frac{e^{-it\xi^2-i\tau(\eta^2-(\eta-\xi)^2-\xi^2)}}{\tau}\,|\xi|^k\,e^{ix\xi}\hat{u_+}(\eta) \widehat{\overline{u}_+}(\eta-\xi)\,d\eta\,d\xi \,d\tau$$
$$=\int\int\frac{e^{-it(\xi^2+2\xi(\eta-\xi))}}{t}\,|\xi|^k\,e^{ix\xi}\,\frac{\hat{u_+}(\eta) \widehat{\overline{u}_+}(\eta-\xi)}{2i\xi(\eta-\xi)}\,d\eta\,d\xi$$
$$+\int_t^{\infty}\int\int\frac{e^{-it\xi^2-2i\tau\xi(\eta-\xi)}}{\tau^2}\,|\xi|^k\,e^{ix\xi}\,\frac{\hat{u_+}(\eta) \widehat{\overline{u}_+}(\eta-\xi)}{2i\xi(\eta-\xi)}\,d\eta\,d\xi\,d\tau.$$
For having a good bound, we need to avoid the powers of $\xi$ in the denominator. That is why the $L^2$ norm was considered apart before this Lemma, and the decay in time obtained was weaker. So by taking $k\geq 1$, we get 
$$\left\|\nabla^ka\int_t^{\infty} e^{i(t-\tau)\partial_x^2}\frac{\left|e^{i\tau\partial_x^2}\omega\right|^2}{\tau}d\tau\right\|_2\lesssim \frac{a}{t}\left\||\xi|^{k-1}\left(\hat{u_+}*\frac{\hat u_+}{\cdot}\right)(\xi)\right\|_{L^2}=\frac{a}{t}\left\|\hat{u_+}\,\mathcal{F}\left(\frac{\hat u_+}{\cdot}\right)\right\|_{\dot{H}^{k-1}}.
$$
For $k=1$ we use the imbedding $L^\infty(\mathbb R)\subset H^1(\mathbb R)$ and we get the upper-bound $\frac{a}{t}\| u_+\|_{H^1}\|u_+\|_{\dot{H}^{-1}}$. For $k\geq 2$ we use the fact that $H^{k-1}(\mathbb R)$ is an algebra, and we get the upper-bound $\frac{a}{t}\|u_+\|_{H^{k-1}}\left(\|u_+\|_{\dot{H}^{-1}}+\|u_+\|_{\dot{H}^{k-2}}\right)$. 

In conclusion, if
$$u_+\in \dot{H}^{-1}\cap \dot{H}^{s-1}.$$
then we can control in $L^2$ the first and higher derivatives up to $s$ of the quadratic terms, like stated in the Lemma.
\end{proof}

The linear term $J_1$ will be treated similarly. We have the following Lemma. 
\begin{lemma}\label{linearosc}
If $u_+\in \dot{H}^{-2}\cap \dot{H}^{s-2}$, then
\begin{equation}\label{linearterm}
\Sigma_{0\leq k\leq s}\left\|\nabla^kJ_1(t)\right\|_2\lesssim \frac{a^2}{t}\left(\left\|u_+\right\|_{\dot{H}^{-2}}+\|u_+\|_{\dot{H}^{s-2}}\right).
\end{equation}
\end{lemma}

\begin{proof}
We write
$$J_1(t)=\int_t^{\infty} e^{i(t-2\tau)\partial_x^2}\left(\frac{\overline{u_+}}{\tau^{1\pm ia^2}}\right)d\tau=\int_t^{\infty} \int\frac{e^{-i(t-2\tau)\xi^2}}{\tau^{1\pm ia^2}}\,e^{ix\xi}\hat{u_+}(-\xi)d\xi d\tau,$$
and we perform an integration by parts in $\tau$,
$$J_1(t)=\int\frac{e^{it\xi^2}}{t^{1\pm ia^2}} \,e^{ix\xi}\frac{\hat{u_+}(-\xi)}{2i\xi^2}d\xi+(1\pm ia^2)\int_t^{\infty}\int\frac{e^{-i(t-2\tau)\xi^2}}{\tau^{2\pm ia^2}} \,e^{ix\xi}\frac{\hat{u_+}(-\xi)}{2i\xi^2}d\xi d\tau.$$
By Plancherel we get
\begin{equation}\label{linear}
\left\|J_1(t)\right\|_2\lesssim \frac{a^2}{t}\left\|\frac{\hat{u_+}(-\xi)}{\xi^2}\right\|_2.
\end{equation}
Since derivatives in the space variable commutes with $e^{i(t-\tau)\partial_x^2}$, we obtain similarly for $k\in \mathbb{N}^*$, 
\begin{equation}\label{linearHs}
\left\|\nabla^kJ_1(t)\right\|_2\lesssim \frac{a^2}{t}\left\|\frac{\hat{u_+}(-\xi)}{|\xi|^{2-k}}\right\|_2,
\end{equation}
and the lemma follows.
\end{proof}

Summarizing, we have obtained that
$$\|A(v)\|_Y\leq C\,\|v\|_Y\left(a^2+\frac{\|u_+\|_{L^1}^2}{t_0}+\frac{(C\,R)^2}{t_0^{\mu+\nu}}\right)\left(1+\underset{t_0\leq t}{\sup}\,\frac{\,t^\mu}{t^{\frac{1}{2}+\nu}}\right)$$
$$+C(u_+)\,\underset{t_0\leq t}{\sup}\,\frac{t^\mu+t^\nu}{t}+C(u_+)\,a\,\,\underset{t_0\leq t}{\sup}\,\frac{t^\nu}{t^\frac 12}+C(u_+)\,a\,\,\underset{t_0\leq t}{\sup}\,\frac{t^\mu}{t}+C(u_+)\,a^2\,\,\underset{t_0\leq t}{\sup}\,\frac{t^\mu+t^\nu}{t},$$
with the constants depending on the $\dot{H}^{-2}\cap H^{s}\cap W^{s,1}$ norm of $u_+$. The choice $\nu=\frac 12$ and $\mu=1$ from the statement of Theorem \ref{longrangeH} satisfy $0<\nu\leq\frac 12$, $0<\mu\leq 1$ and $\mu\leq\frac{1}{2}+\nu$. Therefore there exists a positive constant $R$, small with respect to $t_0$, and a positive constant $a_0$, such that for all $a\leq a_0$, and all $u_+\in \dot{H}^{-2}\cap H^{s}\cap W^{s,1}$ small with respect to $a_0$, $t_0$ and $R$, we can apply the fixed point theorem in $Y_R$. We get this way a solution $v$ of equation \eqref{GP} such that 
$$v-v_1\in \mathcal{C}([t_0,\infty),H^s(\mathbb{R})),$$ 
satisfying 
\begin{equation}\label{finaldecay}
\|(v-v_1)(t)\|_{L^2}\leq\frac{R}{\sqrt{t}}\,\,\,,\,\,\,\|\nabla^k(v-v_1)(t)\|_{L^2}\leq\frac{R}{t},
\end{equation}
and the proof of Theorem \ref{longrangeH} is complete.

\section{The construction of the binormal flow for positive times}\label{flowsection}
\bigskip
In this section we show how to construct a binormal flow for all times $t\geq 0$ from a solution obtained in Theorem \ref{longrangeH}. The next section will contain the proof of the fact that this new binormal flow is close to the selfsimilar one $\chi_a$, and that a singularity is still formed at time $t=0$.

We fix $\epsilon>0$, $a<a_0$ and $t_0=1/\tilde{t_0}$. With the notations of the previous section, let $R$ small enough with respect to $\epsilon$, $t_0$ and $a$, and $u_+\in \dot{H}^{-2}\cap H^{s}\cap W^{s,1}$ small with respect to $\epsilon$, $a$, $t_0$ and $R$. We consider the corresponding  solution $v$ of Theorem \ref{longrangeH} with $s=3$ and equation \eqref{vortex} with focusing sign $+$. We take $s=3$ to have $v$ smooth: $H^3(\mathbb{R})\subset\mathcal{C}^\frac 52(\mathbb{R})$ is enough for our purposes. In subsection \S\ref{estvsection} we give some estimates in time on $v$, $v-a$ and $v-v_1$, obtained from the statement of Theorem \ref{longrangeH}. Then, we define in \S\ref{ctau} a curvature and a torsion, and use these estimates for calculating their leading terms as $t$ goes to zero. These estimates will be used throughout the rest of the paper. The last subsection \S\ref{constrsection} concerns the new binormal flow. Using the curvature and the torsion, that are continuous for $t>0$ but not at $t=0$, the binormal flow is constructed for all positive times. The estimates in \S\ref{ctau} allows us to obtain a limit at $t=0$ for the flow of curves.

\subsection{Estimates from Theorem \ref{longrangeH}}\label{estvsection}
We define $f$ by
$$v(t,y)=a+f(t,y),$$
that is
$$f(t,y)=e^{i a^2\log t}e^{it\partial_x^2}{u_+}(y)+(v-v_1)(t,y).$$
Hereafter in this section, when for a given $h$ we write the expression  $\partial_x h\left(\frac 1t,\frac xt\right)$ we shall mean $ g'(x)$ with $g(x)=h\left(\frac 1t,\frac xt\right)$. 

When $t$ goes to zero, we have different estimates for the two terms of $f\left(\frac 1t,\frac xt\right)$.  For the second one we get from the estimates \eqref{finaldecay} on $v-v_1$,
$$\left\|(v-v_1)\left(\frac 1t,\frac xt\right)\right\|_{L^\infty}\leq R\,t^{\frac 34}\qquad,\qquad\left\|\partial_x (v-v_1)\left(\frac 1t,\frac xt\right)\right\|_{L^\infty}\leq R.$$

For the first term of $f$ we have only the dispersion decay rate 
$$\left\|e^{- i a^2\log t}e^{i\frac 1t\partial_x^2}u_+\left(\frac xt\right)\right\|_{L^\infty}\leq C(u_+)\,\sqrt{t}\qquad,\qquad\left\|\partial_x e^{- i a^2\log t}e^{i\frac 1t\partial_x^2}u_+\left(\frac xt\right)\right\|_{L^\infty}\leq \frac{C(u_+)}{\sqrt{t}}.$$
Since $R$ is small enough with respect to $t_0$ and $a$, and $u_+$ is small with respect to $R$, $t_0$ and $a_0$, we get
\begin{equation}\label{f}
\left\|f\left(\frac 1t,\frac xt\right)\right\|_{L^\infty}\leq C(R)\sqrt{t}
\qquad,\qquad \left\|\partial_x f\left(\frac 1t,\frac xt\right)\right\|_{L^\infty}\leq \frac{C(R)}{\sqrt{t}}.\end{equation}

However, at $x=0$ we get a better decay for the first derivative. From the expression of the free Schr\"odinger evolution,
$$ e^{- i a^2\log t}e^{i\frac 1t\partial_x^2}u_+\left(\frac xt\right)= e^{- i a^2\log t}\frac{e^{i\frac{x^2}{4t}}}{\sqrt{\frac it}}\int e^{-i\frac{xy}{2}}e^{i\frac{y^2}{4}t}u_+(y)dy,$$
and taking $x=0$,
$$\left|\partial_x e^{- i a^2\log t}e^{i\frac 1t\partial_x^2}u_+\left(\frac xt\right)\right|_{x=0}=\sqrt{t}\left|\int \frac y2 e^{i\frac{y^2}{4}t}u_+(y)dy\right|$$
$$\leq \sqrt{t}\|yu_+\|_{L^1}\leq\sqrt{t}\left(\|y^2u_+\|_{L^2}+\|u_+\|_{L^2}\right).$$
Therefore at $x=0$,
\begin{equation}\label{cx}
\left|\partial_x f\left(\frac 1t,\frac xt\right)\right|_{x=0}\leq C(R).
\end{equation}

\subsection{The curvature and the torsion}\label{ctau}
We start by defining a curvature and a torsion from the solution $v$ of Theorem \ref{longrangeH}.
Let us recall that since $v$ is a solution of the focusing equation \eqref{GP}, then its pseudo-conformal transform  
$$u(t,x)=\D\frac{e^{i\frac{x^2}{4t}}}{\sqrt{t}}\overline{v}\left(\frac 1t,\frac xt\right),$$
is a solution of 
$$iu_t+u_{xx}+\left(|u|^2-\frac{a^2}{t}\right)u=0.$$
Since $v$ is regular enough and does not vanish, 
we can define two real functions $\tau$ and $\phi$ such that $u$ 
$$u(t,x)=c(t,x)e^{i\phi(t,x)}.$$
We define $\tau(t,x):=\phi_x(t,x)$, so
$$u(t,x)=c(t,x)e^{i\int_0^x\tau(t,s)ds +\phi(t,0)}.$$

Then, the function
$$\tilde{u}(t,x)=c(t,x)e^{i\int_0^x\tau(t,s)ds},$$
is a filament function, solution of \eqref{Hasimoto-eq} with $A(t)$ replaced by $\frac{a^2}{t}+\phi_t(t,0)$. As will be seen in the next subsection \S\ref{constrsection}, there exists a binormal flow of curves such that the curvature and the torsion are $c$ and $\tau$.\\

As $t$ goes to zero, we shall compute the leading terms of $(c,\tau)$. We have
$$c(t,x)=|u(t,x)|=\frac{1}{\sqrt{t}}\left|v\left(\frac 1t,\frac xt\right)\right|,$$
and
$$\tau(t,x)=\Im\frac{u_x(t,x)}{u(t,x)}=\Im\frac{\frac{ix}{2t}\overline{v}\left(\frac 1t,\frac xt\right)+\partial_x \overline{v}\left(\frac 1t,\frac xt\right)}{\overline{v}\left(\frac 1t,\frac xt\right)}.$$

Since $v=a+f$, the square of the curvature is 
$$c^2(t,x)=\frac {a^2}{t}\left(1+ \frac 2a \Re f+\frac{|f|^2}{a^2}\right)\left(\frac 1t,\frac xt\right).$$
Because $c$ and $a$ are positive, 
$$
c(t,x)-\frac{a}{\sqrt{t}}=\frac{1}{c+a/\sqrt t}\left( \frac {2a}{ t} \Re f+\frac{|f|^2}{t}\right)\left(\frac 1t,\frac xt\right),
$$
and in view of estimate \eqref{f} on $f$, we obtain the estimate on the curvature,
\begin{equation}\label{c}
\left|c(t,x)-\frac{a}{\sqrt{t}}\right|\leq\frac {1}{\sqrt t}\left| \left( 2\Re f+\frac{|f|^2}{a}\right)\left(\frac 1t,\frac xt\right)\right|\leq C(R).
\end{equation}
We recall that $R$ is small with respect to $t_0$ and $a$. 
It follows that $c>a/2\sqrt t$. Hence similarly, from $\partial_x c^2=2c\,\partial_x c$, we get an estimate on the first derivative in space of the curvature, by using \eqref{f},
\begin{equation}\label{c_x}
|\partial_x c|\leq\frac {1}{\sqrt t}\left| \partial_x\left( 2\Re f+\frac{|f|^2}{a}\right)\left(\frac 1t,\frac xt\right)\right|\leq \frac{C(R)}{t}.
\end{equation}
At $x=0$ we can use \eqref{cx} and get
\begin{equation}\label{c_xzero}
|\partial_x \,c (t,0)|\leq\frac {1}{\sqrt t}\left| \partial_x\left( 2\Re f+\frac{|f|^2}{a}\right)\left(\frac 1t,\frac xt\right)\right|\leq \frac{C(R)}{\sqrt{t}}.
\end{equation}

The torsion is well defined and is given by
$$\tau(t,x)=\Im\frac{\frac{ix}{2t}\left(a+\overline{f}\left(\frac 1t,\frac xt\right)\right)+\partial_x \overline{f}\left(\frac 1t,\frac xt\right)}{a+\overline{f}\left(\frac 1t,\frac xt\right)}.$$
Then
$$\tau(t,x)-\Im\left(\frac{ix}{2t}-\frac 1a\partial_x\overline{f}\left(\frac 1t,\frac xt\right)\right)=-\Im\frac{\partial_x\overline{f}\left(\frac 1t,\frac xt\right) \overline{f}\left(\frac 1t,\frac xt\right)}{a\left(a+\overline{f}\left(\frac 1t,\frac xt\right) \right)},$$
and so we get,
\begin{equation}\label{tau}
\left|\tau(t,x)-\frac{x}{2t}-\Im\left(\frac{1}{a}\partial_x f\left(\frac 1t,\frac xt\right)\right)\right|\leq\frac{2}{a^2}\left|\partial_x f\left(\frac 1t,\frac xt\right) f\left(\frac 1t,\frac xt\right)\right|\leq C(R)\sqrt{t}.
\end{equation}
In particular, by \eqref{f} we get
\begin{equation}\label{taugen}
\left|\tau(t,x)-\frac{x}{2t}\right|\leq \frac{C(R)}{\sqrt{t}},
\end{equation}
and by \eqref{cx} we have at $x=0$,
\begin{equation}\label{tauzero}
\left|\tau(t,0)\right|\leq C(R).
\end{equation}

Let us give also an estimate at $x=0$. By the definition of $u$ and $\phi$ we have 
$$c(t,0)e^{i\phi(t,0)}-\frac{a}{\sqrt{t}}=u(t,0)-\frac{a}{\sqrt{t}}=\frac{1}{\sqrt{t}}(v-a)\left(\frac 1t,0\right)=\frac{1}{\sqrt{t}}f\left(\frac 1t,0\right).$$
From estimates \eqref{f} on $f$ we get that
$$\left|c(t,0)e^{i\phi(t,0)}-\frac{a}{\sqrt{t}}\right|\leq C(R),$$
and by using \eqref{c},
$$\frac{a}{\sqrt{t}}\left|e^{i\phi(t,0)}-1\right|\leq \left|\frac{a}{\sqrt{t}}-c(t,0)\right|+\left|c(t,0)e^{i\phi(t,0)}-\frac{a}{\sqrt{t}}\right|\leq C(R).$$
Therefore, since $R$ is small enough, 
\begin{equation}\label{phase}
\left|e^{i\frac{\phi(t,0)}{2}}-1\right|=\frac{\left|e^{i\phi(t,0)}-1\right|}{\left|e^{i\frac{\phi(t,0)}{2}}+1\right|}\leq C(R)\sqrt{t}.\end{equation}

Finally, let us recall that the curvature and the torsion of the selfsimilar binormal flow $\chi_a$ are
\begin{equation}\label{ctaua}
c_a(t,x)=\frac{a}{\sqrt{t}},\qquad ,\qquad \tau_a(t,x)=\frac{x}{2t}.
\end{equation}
Therefore all the estimates in this subsection show that $(c,\tau)$ is uniformly close to $(c_a,\tau_a)$. This will be used in the next section \S\ref{formationsection}.

\subsection{The integration of the binormal flow}\label{constrsection} 
From the curvature and the torsion defined in the previous subsection, we shall construct a corresponding family of curves solution of \eqref{binormal}. 
We first construct its tangent, normal and binormal vectors $(T,n,b)(t,x)$ in the following way. For a given $(T,n,b)(\tilde{t_0},0)$, we define $(T,n,b)(t,0)$ by imposing 
\begin{equation}\label{timeFrenet}
\left(\begin{array}{c}
T\\ n\\ b
\end{array}\right)_t(t,0)=
\left(\begin{array}{ccc}
0 & -c\,\tau & c_x \\ c\,\tau & 0 & \left(\frac{c_{xx}-c\tau^2}{c}\right) \\ -c_x &  -\left(\frac{c_{xx}-c\tau^2}{c}\right) & 0 
\end{array}\right)
\left(\begin{array}{c}
T\\ n\\ b
\end{array}\right)(t,0).
\end{equation}
This is the system that the time derivatives of the tangent, normal, and binormal of a binormal flow verifies. This will be proved in the Appendix. 

Then, we construct $(T,n,b)(t,x)$ from $(T,n,b)(t,0)$ by integrating the Frenet system for fixed $t$
$$\left(\begin{array}{c}
T\\ n\\ b
\end{array}\right)_x(t,x)=
\left(\begin{array}{ccc}
0 & c & 0 \\ -c & 0 & \tau \\ 0 &  -\tau & 0 
\end{array}\right)
\left(\begin{array}{c}
T\\ n\\ b
\end{array}\right)(t,x).
$$

 This way $T$ will solve, (see the Appendix)
$$T_t=T\land T_{xx}.$$

With $T$ constructed this way, for a given curve $\chi(\tilde{t_0},0)$, we define for all $\tilde{t_0}>t>0$,
$$\chi(t,x):=\chi(\tilde{t_0},0)-\int^{\tilde{t_0}}_t (cb)(t',0)dt'+\int_0^xT(t,s)ds.$$
Using the Frenet system,
$$T_t=T\land T_{xx}=T\land (cn)_x=T\land(c_x n +c\tau b)=-c\tau n+c_xb,$$
and it follows that $\chi$ solves the binormal flow equation \eqref{binormal1}.\\

Therefore $\chi(t,x)$ is constructed for all times when the curvature and the torsion are regular, that is for $t>0$. Finally, by using \eqref{binormal1} and the expression of the curvature \eqref{c}  we have
\begin{equation}\label{tzero}
\left|\chi(t_1,x)-\chi(t_2,x)\right|=
\left|\int_{t_1}^{t_2} c(t,x)b(t,x) dx\right|\leq \int_{t_1}^{t_2} \frac{Ca}{\sqrt{t}} dt\underset{t_1,t_2\rightarrow 0}{\longrightarrow} 0.
\end{equation}
By denoting $\chi_0(x)$ the limit at $t=0$, we obtain similarly that for all $x\in(-\infty,\infty)$,
$$\left|\chi(t,x)-\chi_0(x)\right|\leq Ca\sqrt{t},$$
and the first part of Theorem \ref{stab} is proved.

\section{The formation of the singularity for the binormal flow}\label{formationsection}
In this section we shall prove the second part of the statement of Theorem \ref{stab}. We shall show that the binormal flow $\chi$ constructed in the previous section is close to the selfsimilar one $\chi_a$. This will allow us to conclude that a corner is still formed at time zero at $x=0$. 

To this purpose, we start by showing that the tangent $T$ of $\chi$ remains close to $T_a$, the tangent of $\chi_a$. This will be done in three steps in the next three subsections. First we show that $(T,n,b)(t,x)$ remains close to $(T_a,n_a,b_a)(t,x)$ at $t=0$. Using this we show in the second step that $(T,n,b)(t,x)$ remains close to $(T_a,n_a,b_a)(t,x)$ for $x\lesssim\sqrt{t}$. In particular, $T(t,x)$ is close to $T_a(t,x)$ for $x\lesssim\sqrt{t}$. This will imply in the final step that $T(t,x)$ is close to $T_a(t,x)$ also for $\sqrt{t}\lesssim x$.

In the last subsection the information that $T(t,x)$ is close to $T_a(t,x)$ is used to show that $\chi(0,x)$ is close to $\chi_a(0,x)$.

\subsection{Estimates at $(t,0)$} Let us recall that in subsection \ref{constrsection} we have constructed $(T,n,b)(t,0)$ by imposing \eqref{timeFrenet},
$$\left(\begin{array}{c}
T\\n\\b
\end{array}\right)_t(t,x)=
\left(\begin{array}{ccc}
0 & -c\,\tau & c_x \\ c\,\tau & 0 &  \frac{c_{xx}-c\tau^2}{c}
\\ -c_x &  -\frac{c_{xx}-c\tau^2}{c} & 0 
\end{array}\right)
\left(\begin{array}{c}
T\\ n\\ b
\end{array}\right)(t,x).$$
As notice when the curvature and the torsion have been defined in subsection \ref{ctau}, the function
$$\tilde{u}(t,x)=c(t,x)e^{i\int_0^x\tau(t,s)ds},$$
is a filament function, solution of \eqref{Hasimoto-eq} with $A(t)$ replaced by $\frac{a^2}{t}+\phi_t(t,0)$. It follows that we have the condition \eqref{Hasimoto-coeff}
$$\frac{a^2}{t}+\partial_t\phi(t,0)=2\left(\frac{c_{xx}-c\tau^2}{c}\right)(t,0)+c(t,0)^2.$$
Therefore we obtain
$$\left(\begin{array}{c}
T\\n\\b
\end{array}\right)_t(t,0)=
\left(\begin{array}{ccc}
0 & -c\,\tau & c_x \\ c\,\tau & 0 &  \frac{c_a^2-c^2}{2}+\frac{\phi_{t}}{2}
\\ -c_x &  -\frac{c_a^2-c^2}{2}-\frac{\phi_{t}}{2} & 0 
\end{array}\right)
\left(\begin{array}{c}
T\\ n\\ b
\end{array}\right)(t,0).$$
In order to get rid of the term $\phi_t(t,0)$, we introduce
$$\tilde{n}+i\tilde{b}=e^{i\frac{\phi}{2}}(n+ib).$$
A straightforward computation gives us
$$\tilde{n}_t+i\tilde{b}_t=e^{i\frac{\phi}{2}}\left(n_{t}+ib_{t}+i\frac{\phi_{t}}{2}n-
\frac{\phi_{t}}{2}b\right)=e^{i\frac{\phi}{2}}(c_{\tau}-ic_{x})T-i\frac{c_a^2-c^2}{2}(\tilde{n}+i\tilde{b}),$$
and
$$\left(\begin{array}{c}
T\\\tilde{n}\\\tilde{b}
\end{array}\right)_t(t,0)=
\left(\begin{array}{ccc}
0 & -c\,\tau\cos\frac{\phi}{2}-c_x\sin\frac{\phi}{2} & -c\,\tau\sin\frac{\phi}{2}+c_x\cos\frac{\phi}{2} 
\\ c\,\tau\cos\frac{\phi}{2}+c_x\sin\frac{\phi}{2} & 0 &  \frac{c_a^2-c^2}{2}
\\ c\,\tau\sin\frac{\phi}{2}-c_x\cos\frac{\phi}{2}  &  -\frac{c_a^2-c^2}{2}& 0 
\end{array}\right)
\left(\begin{array}{c}
T\\ \tilde{n}\\\tilde{b}
\end{array}\right)(t,0).$$

We choose as an initial data $(T,\tilde{n},\tilde{b})(\tilde{t}_0,0)=(T_a,n_a,b_a)(\tilde{t}_0,0)$. Since $(T_a,n_a,b_a)(t,0)$ is the orthonormal basis of $\mathbb R^3$, we obtain  
$$\left|\begin{array}{c}
T-T_a\\\tilde{n}-n_a\\\tilde{b}-b_a
\end{array}\right|(t,0)\leq 3\int_t^{\tilde{t_0}} \left(|c\,\tau|+|c_x|+|c_a^2-c^2|\right)(\sigma,0)\,d\sigma.$$
From the expressions \eqref{c},\eqref{c_xzero},\eqref{tauzero} of the curvature and the torsion at $x=0$, 
$$\left(|c\,\tau|+|c_x|+|c_a^2-c^2|\right)(\sigma,0)\leq \frac{C(R)}{\sqrt{\sigma}} .$$
Therefore we get
\begin{equation}\label{prestep1}
\left|\begin{array}{c}
T-T_a\\\tilde{n}-n_a\\\tilde{b}-b_a
\end{array}\right|(t,0)\leq C(R),
\end{equation}
and the fact that $(T,\tilde{n},\tilde{b})(t,0)$ has a limit as $t$ goes to zero. Finally,
$$|(n+ib)-(n_a+ib_a)|\leq|(n+ib)-e^{i\frac{\phi}{2}}(n+ib)|+|(\tilde{n}+i\tilde{b})-(n_a+ib_a)|,$$
and in view of \eqref{phase} and \eqref{prestep1}, we obtain
\begin{equation}\label{step1}
\left|\begin{array}{c}
T-T_a\\n-n_a\\b-b_a
\end{array}\right|(t,0)\leq C(R),
\end{equation}
and the fact that $(T,n,b)(t,0)$ has a limit as $t$ goes to zero.

\subsection{Estimates at $(t,x)$ for $x\lesssim\sqrt{t}$}

Let us denote 
$$N=n+ib.$$ 
The tangent, normal and binormal were constructed in subsection \ref{constrsection} such that we can use the Frenet system. 
This gives us
$$N_x=-cT+\tau b-i\tau n=-cT-i\tau N,$$
so that
$$(N-N_a)_x=-(c-c_a)T-c_a(T-T_a)-i(\tau-\tau_a)N-i\tau_a(N-N_a).$$
In particular,
\begin{equation}\label{N}
e^{-i\frac{x^2}{4t}}\left(e^{i\frac{x^2}{4t}}(N-N_a)\right)_x=-(c-c_a)T-c_a(T-T_a)-i(\tau-\tau_a)N,
\end{equation}
and
$$(T-T_a)_x=cn-c_an_a=(c-c_a)n+c_a(n-n_a).$$
If we denote
$$\Sigma^2=|T-T_a|^2+|N-N_a|^2,$$
we can compute using \eqref{N}
$$\Sigma^2_x=2<(c-c_a)n+c_a(n-n_a),T-T_a>$$
$$+2<-(c-c_a)T-c_a(T-T_a)-i(\tau-\tau_a)N,N-N_a>.$$
The tangent and the normal vectors are of norm $1$, and $|N|$ is bounded  by $2$, so
$$2\Sigma\Sigma_x\leq2|c-c_a||T-T_a|+c_a|n-n_a||T-T_a|$$
$$+2|c-c_a||N-N_a|+c_a|T-T_a||N-N_a|+|\tau-\tau_a||N-N_a|$$
$$\leq 2(|c-c_a|+|\tau-\tau_a|)\Sigma+c_a\Sigma^2.$$
Therefore 
$$\left(e^{-x\frac{a}{2\sqrt{t}}}\Sigma\right)_x\leq e^{-x\frac{a}{2\sqrt{t}}}(|c-c_a|+|\tau-\tau_a|),$$
and so
$$\Sigma(t,x)\leq e^{x\frac{a}{2\sqrt{t}}}\Sigma(t,0)+e^{x\frac{a}{2\sqrt{t}}}\int_0^x e^{-y\frac{a}{2\sqrt{t}}}(|c-c_a|+|\tau-\tau_a|)(t,y)dy.$$
For $x<M\sqrt{t}$, with $M$ to be chosen later, we have
$$\Sigma(t,x)\leq e^{\frac{Ma}{2}}\Sigma(t,0)+e^{\frac{Ma}{2}}M\sqrt{t}\underset{0\leq y\leq M\sqrt{t}}{\sup}(|c-c_a|+|\tau-\tau_a|)(t,y).$$
Using the bounds \eqref{c}, \eqref{taugen} of the curvature and torsion, 
$$\Sigma(t,x)\leq e^{\frac{Ma}{2}}\Sigma(t,0)+e^{\frac{Ma}{2}}M\,C(R).$$
In the last subsection we shall choose $M$ large, so we can write
$$\Sigma(t,x)\leq e^{Ma}\Sigma(t,0)+e^{Ma}C(R).$$
By  \eqref{step1}, that is proved in the previous subsection, $\Sigma(t,0)\leq C(R)$ and we get
$$\Sigma(t,x)\leq e^{Ma}C(R).$$
In conclusion, in the region $x\leq M\sqrt{t}$
\begin{equation}\label{step2}
\left|\left(\begin{array}{c}
T-T_a\\n-n_a\\b-b_a
\end{array}\right)(t,x)\right|\leq e^{Ma}C(R).
\end{equation}

\subsection{Estimates at $(t,x)$ for $\sqrt{t}\lesssim x$ }\label{step3section}
Using the Frenet system again we write
$$T(t,x')-T(t,x)=\int_{x}^{x'}cn=\int_{x}^{x'}c\left(1-\frac{\tau}{\tau_a}\right)n+\frac{c}{\tau_a}\tau n.$$
Since $b_x=-\tau n$, we do an integration by parts in the last term,
$$T(t,x')-T(t,x)=\left[-\frac{c}{\tau_a}b\right]_{x}^{x'}+\int_{x}^{x'}c\frac{\tau_a-\tau}{\tau_a}n+\left(\frac{c}{\tau_a}\right)'b,$$
and by using the explicit expression $\tau_a(t,x)=\frac{x}{t}$,
\begin{equation}\label{idT}
T(t,x')=T(t,x)+\frac{2t\,c(t,x)}{x}b(t,x)-\frac{2t\,c(t,x')}{x'}b(t,x')
\end{equation}
$$+\int_{x}^{x'} \frac{2t\, c}{s}\left(\frac{s}{2t}-\tau\right)n+\left(\frac{2tc}{s}\right)_sb\,ds.$$
Now we write the difference
$$(T-T_a)(t,x')=(T-T_a)(t,x)+\frac{2t(c-c_a)}{x}b+\frac{2tc_a}{x}(b-b_a)-\frac{2t(c-c_a)}{x'}b-\frac{2tc_a}{x'}(b-b_a)$$
$$+\int_x^{x'} \frac{2t\, c}{s}\left(\frac{s}{2t}-\tau\right)n+\left(\frac{2t(c-c_a)}{s}\right)_sb+\left(\frac{2tc_a}{s}\right)_s(b-b_a)\,ds.$$
By choosing $x'>x=M\sqrt{t}$,
$$|T-T_a|(t,x')\leq |T-T_a|(t,M\sqrt{t})+\frac{4\sqrt{t}|c-c_a|}{M}+\frac{8\sqrt{t}c_a}{M}$$
$$+\left|\int_{M\sqrt{t}}^{x'} \frac{2t\, c}{s}\left(\frac{s}{2t}-\tau\right)n+\left(\frac{2t(c-c_a)}{s}\right)_sb+\left(\frac{2tc_a}{s}\right)_s(b-b_a)\,ds\right|.$$
The result \eqref{step2} of the previous subsection together with the decay \eqref{c} of $c-c_a$ give us
\begin{equation}\label{diffT}
|(T-T_a)(t,x')|\leq e^{Ma}\,C(R)+\frac{8a}{M}(1+C(R))
\end{equation}
$$+\left|\int_{M\sqrt{t}}^{x'} \frac{2t\, c}{s}\left(\frac{s}{2t}-\tau\right)n+\left(\frac{2t(c-c_a)}{s}\right)_sb+\left(\frac{2tc_a}{s}\right)_s(b-b_a)\,ds\right|.$$
We denote $I_1$, $I_2$ and $I_3$ the integral terms.

The last term can be easily estimated by
$$|I_3|=\left|\int_{M\sqrt{t}}^{x'}\left(\frac{2tc_a}{s}\right)_s(b-b_a)\,ds\right|\leq 2\left|\int_{M\sqrt{t}}^{x'}\left(\frac{2a\sqrt{t}}{s}\right)_s\right|<\frac{4a}{M}.$$
Now we consider the second term in the integral
$$|I_2|=\left|\int_{M\sqrt{t}}^{x'}\left(\frac{2t(c-c_a)}{s}\right)_sb\,ds\right|\leq \int_{M\sqrt{t}}^{x'}\left|\frac{2t\,c_s}{s}\right|+\left|\frac{2t\,(c-c_a)}{s^2}\right|\,ds.$$
We have from the estimates \eqref{c} and \eqref{c_x} on the $c-c_a$ and on $c_x$ respectively,
$$\left|\frac{2t\,c_s(t,s)}{s}\right|\leq\frac{C\sqrt{t}}{s}\left|\partial_sf\left(\frac 1t,\frac st\right)\right|+\frac{C\sqrt{t}}{s}\left|\partial_s|f|^2\left(\frac 1t,\frac st\right)\right|.$$
and
$$\left|\frac{2t\,(c-c_a)(t,s)}{s^2}\right|\leq\frac{2\sqrt{t}}{s^2}|f|\left(\frac 1t,\frac st\right)+\frac{\sqrt{t}|f|^2\left(\frac 1t,\frac st\right)}{s^2}.$$
By using Cauchy-Schwarz' inequality and the bounds \eqref{f} on $f$, we get
$$\left|I_2\right|\leq  \int_{M\sqrt{t}}^{x'}\frac{C\sqrt{t}}{s}\left|\partial_sf\left(\frac 1t,\frac st\right)\right|\,ds+\frac{C(R)}{M}.$$
Finally from the expression \eqref{tau} of the torsion,
$$|I_1|=\left|\int_{M\sqrt{t}}^{x'}\frac{2tc(t,s)}{s}\left(\frac{s}{2t}-\tau(t,s)\right)\right|\leq C\left|\int_{M\sqrt{t}}^{x'}\frac{\sqrt t}{s}\partial_s f\left(\frac 1t,\frac st\right)\left(1+\frac1a f\left(\frac 1t,\frac st\right)\right)\right|$$
$$\leq  \int_{M\sqrt{t}}^{x'}\frac{C\sqrt{t}}{s}\left|\partial_sf\left(\frac 1t,\frac st\right)\right|\,ds+\frac{C(R)}{M}.$$

In conclusion, we can transform \eqref{diffT} into
\begin{equation}\label{diffT2}
|(T-T_a)(t,x')|\leq e^{Ma}\mathcal{O}(R,u_+)+\frac{12a}{M}(1+\mathcal{O}(R,u_+))+ \int_{M\sqrt{t}}^{x'}\frac{C\sqrt{t}}{s}\left|\partial_sf\left(\frac 1t,\frac st\right)\right|\,ds.
\end{equation}

Next we need the following lemma.

\begin{lemma}\label{eststep3}
The following estimate holds
$$\int_{x}^{x'}\frac{C\sqrt{t}}{s}\left|\partial_sf\left(\frac 1t,\frac st\right)\right|\,ds\leq \|\hat{u}_+\|_{L^1(x,x')}+\frac {C(R)}{\sqrt{x}}\,t.$$
\end{lemma}

\begin{proof}
Recall that $f$ is given by
$$f(t,y)=e^{i a^2\log t}e^{it\partial_x^2}{u_+}(y)+(v-v_1)(t,y).$$
Denote by $A_1$ and $A_2$ the  two terms in the above sum.

The second one can be estimated easily by using Cauchy-Schwarz inequality,
$$A_2=\int_{x}^{x'}\frac{\sqrt{t}}{s}\left|\partial_s(v-v_1)\left(\frac 1t,\frac st\right)\right|\,ds\leq \sqrt{t}\left\|\frac{1}{s}\right\|_{L^2(x,x')}\left\|\partial_s(v-v_1)\left(\frac 1t,\frac st\right)\right\|_2,$$
and the rate of decay of Theorem \ref{longrangeH},
$$A_2\leq \sqrt{t}\frac{1}{\sqrt{x}}\frac{1}{t}\sqrt{t}\,C(R)\,t=\frac{C(R)}{\sqrt{x}}\,t.$$

By using the expression of the free Schr\"odinger evolution,
$$A_1=\int_{x}^{x'}\frac{\sqrt{t}}{s}\left|\partial_s\left( \frac{e^{i\frac{s^2}{4t}}}{\sqrt{\frac it}}\int e^{-i\frac{sy}{2}}e^{i\frac{y^2}{4}t}u_+(y)dy\right)\right|ds=\int_{x}^{x'}\frac{t}{s}\left|\partial_s\left(e^{i\frac{s^2}{4t}}\int e^{-i\frac{sy}{2}}e^{i\frac{y^2}{4}t}u_+(y)dy\right)\right|ds$$
$$\leq\frac 12\int_{x}^{x'}\left|\int e^{-i\frac{sy}{2}}e^{i\frac{y^2}{4}t}u_+(y)dy\right|ds+\int_{x}^{x'}\frac{t}{s}\left|\int e^{-i\frac{sy}{2}}e^{i\frac{y^2}{4}t}y u_+(y)dy\right|ds$$
$$\leq \int_{x}^{x'}\left|\hat u_+\left(\frac s2\right)\right|\frac{ds}{2}+\int_{x}^{x'}\left|\int e^{-i\frac{sy}{2}}\left(e^{i\frac{y^2}{4}t}-1\right)u_+(y)dy\right|\frac{ds}{2}+t\left\|\frac{1}{s}\right\|_{L^2(x,x')}\| y u_+\|_{L^2}.$$
In the second term we perform an integration by parts, and we obtain
$$A_1\leq  \|\hat{u}_+\|_{L^1(x,x')}+\int_{x}^{x'}\frac 1s\left|\int e^{-i\frac{sy}{2}}\partial_y\left(\left(e^{i\frac{y^2}{4}t}-1\right)u_+(y)\right)dy\right|ds+\frac{t}{\sqrt{x}}\| y u_+\|_{L^2}$$
$$\leq \|\hat{u}_+\|_{L^1(x,x')}+\left\|\frac{1}{s}\right\|_{L^2(x,x')}\left\|\partial_y\left(\left(e^{i\frac{y^2}{4}t}-1\right)u_+(y)\right) \right\|_{L^2}+\frac{t}{\sqrt{x}}\| y u_+\|_{L^2}$$
$$\leq \|\hat{u}_+\|_{L^1(x,x')}+\frac{1}{\sqrt{x}}\left\|\left(e^{i\frac{y^2}{4}t}-1\right)\partial_yu_+(y)\right\|_{L^2}+\frac{t}{\sqrt{x}}\|yu_+\|_{L^2}+\frac{t}{\sqrt{x}}\| y u_+\|_{L^2}$$
Since $(1+y^2)u_+$ is in $H^1$, the lemma follows.   
\end{proof}

Let us first notice that in 1-D we can upper bound $\|\hat{u}_+\|_{L^1}\leq \|u_+\|_{\dot{H}^1}$. Since $M$ will be chosen large, the Lemma allows us to re-write \eqref{diffT2} for all $x'\geq x=M\sqrt{t}$, 
$$|(T-T_a)(t,x')|\leq \frac {12a}{M}+e^{Ma}\,C(R).$$

\subsection{The formation of the singularity}\label{cornersection}
Putting together the results of the three previous subsection, we have obtained that for all $\epsilon>0$, and choosing first $M$ large in terms of $\epsilon$ and $a$, then $R$, and then $u_+$ small in terms of $\epsilon$, $a$ and $t_0$, we get that for all $x$ and as $t$ goes to zero,
$$|T(t,x)-T_a(t,x)|\leq\epsilon.$$

Also notice that  Lipschitz property of $\chi_0$ easily follows from \eqref{tzero}.

For $t>0$, $x>0$, 
$$\chi(t,x)-\chi(t,0)=\int_0^x T(t,s)\,ds$$
$$= A^+_ax +\int_0^x (T(t,s)-T_a(t,s))\,ds+\int_0^x (T_a(t,s)-A^+_a)\,ds,$$
so that,
$$|\chi(t,x)-\chi(t,0)-A^+_ax-\int_0^x (T_a(t,s)-A^+_a)\,ds|\leq\epsilon x.$$
As it was said in the Introduction the behaviour of $T_a(t,s)$ as $t$ goes to zero was studied in \cite{SJL}. We shall use \eqref{ii} with $\psi=\mathbb{I}_{[0,x]}$ to get 
$$\lim_{t\rightarrow 0^+}\int_0^x (T_a(t,s)-A^+_a)\,ds=0.$$
Therefore, by letting $t$ going to zero we get 
$$|\chi_0(x)-\chi_0(0)-A^+_ax|\leq\epsilon\,x,$$
and the proof of Theorem \ref{longrangeH} is complete.

\section{Further properties of the binormal flow}\label{extrasection}
In this section we shall prove that the tangent vector $T(t,x)$ has a limit, for fixed $t$, as $|x|$ goes to infinity.

We have obtained in the previous subsection the identity \eqref{idT},
$$T(t,x')-T(t,x)=\frac{2t\,c(t,x)}{x}b(t,x)-\frac{2t\,c(t,x')}{x'}b(t,x')$$
$$+\int_{x}^{x'} \frac{2t\, c}{s}\left(\frac{s}{2t}-\tau\right)n-\left(\frac{2tc}{s}\right)_sb\,ds.$$
We shall prove that the difference $T(t,x')-T(t,x)$ goes to zero as $x,x'$ go to infinity. We recall the expression \eqref{c} of the curvature,
$$2t\,c(t,s)\leq 2\sqrt{t}a+2\sqrt{t}|f|\left(\frac 1t,\frac st\right)+2\sqrt{t}\frac{|f|^2\left(\frac 1t,\frac st\right)}{a}.$$
In view of the estimates \eqref{f} on $f$, when $x,x'$ go to infinity, 
$$|T(t,x')-T(t,x)|\lesssim\int_{x}^{x'} \left|\frac{2t\, c}{s}\left(\frac{s}{2t}-\tau\right)\right|+\left|\frac{2t\,c_s(t,s)}{s}\right|+\left|\frac{2t\,c(t,s)}{s^2}\right|\,ds=B_1+B_2+B_3$$

As before, the last term $B_3$ is integrable on $(x,x')$, and its integral goes to zero as $x,x'$ go to infinity. 

By using the expression \eqref{cx} of the derivative of the curvature, 
$$B_2=\int_{x}^{x'}\left|\frac{2t\,c_s(t,s)}{s}\right|\leq\int_{x}^{x'}\frac{C\sqrt{t}}{s}\left|\partial_sf\left(\frac 1t,\frac st\right)\right|+\int_{x}^{x'}\frac{C\sqrt{t}}{s}\left|\partial_s|f|^2\left(\frac 1t,\frac st\right)\right|.$$
We apply Cauchy-Schwarz inequality in both integrals, and use the estimates \eqref{f} on $f$. On one hand,
$$\int_x^{x'}\left|\frac{2\sqrt{t}}{s}\partial_sf\left(\frac 1t,\frac st\right)\right|dy\leq \sqrt{t}\left\|\frac{1}{s}\right\|_{L^2(x,x')}\left\|\partial_sf\left(\frac 1t,\frac st\right)\right\|_{L^2}\leq \left\|\frac{1}{s}\right\|_{L^2(x,x')}C(R),$$
so as $x,x'$ go to infinity, this integral goes to zero. On the other hand,
$$\int_x^{x'}\frac{C\sqrt{t}}{s}\left|\partial_s|f|^2\left(\frac 1t,\frac st\right)\right|\leq C\frac{\sqrt{t}}{x}\left\|f\left(\frac 1t,\frac st\right)\right\|_{L^2}\left\|\partial_sf\left(\frac 1t,\frac st\right)\right\|_{L^2}\leq \frac{\sqrt{t}}{x}C(R),$$
so $B_2$ goes to zero as $x,x'$ go to infinity. 

Finally, by using the torsion expression \eqref{tau}, $$B_1=\int_{x}^{x'}\left|\frac{2tc(t,s)}{s}\left(\frac{s}{2t}-\tau(t,s)\right)\right|\leq C\left|\frac{\sqrt t}{s}\partial_s f\left(\frac 1t,\frac st\right)\left(1+\frac1a f\left(\frac 1t,\frac st\right)\right)\right|,$$
so we can treated $B_1$ similarly. 

In conclusion, the difference $T(t,x')-T(t,x)$ goes to zero as $x',x$ go to infinity. It follows that for all times there is a limit
$$A^+(t)= \underset{x\rightarrow\infty}{\lim}T(t,x).$$
The same argument can be done as $x$ goes to $-\infty$.

\section{Appendix}
We recall here some general facts about the binormal flow. We show how from a curvature and a torsion defined from a filament function, one can construct tangent, normal and binormal vectors with the properties required by a binormal flow, and in particular \eqref{map}, 
$$T_t=T\land T_{xx}.$$

First, we shall compute the system of the derivatives in time of $(T,n,b)$, the tangent, normal and binormal vectors of a general binormal flow of curves.  Since we have \eqref{map}, by using the Frenet system it follows that
$$T_t=T\land T_{xx}=T\land (cn)_x=T\land(c_x n +c\tau b)=- c\tau n+c_x b.$$
From
$$T_x=cn$$
we get
$$cn_t=-c_tn+(T_t)_x=-c_tn+(c_xb-c\tau n)_x$$
$$=-c_tn+c_{xx}b-c_x\tau n-(c\tau)_xn+c^2\tau T-c\tau^2 b.$$
The vector $n$ is unitary, so $<n_t,n>=0$. Hence $n_t$ is decomposed only in $T$ and $b$. We have
$$n_t=c\tau T+\left(\frac{c_{xx}-c\tau^2}{c}\right)b.$$
Therefore, since $(T,n,b)$ form an orthonormal basis of $\R^3$, the system of derivatives in time of the  tangent, normal and binormal vectors of a binormal flow is
\begin{equation}\label{Frenettimebis}
\left(\begin{array}{c}
T\\ n\\ b
\end{array}\right)_t(t,x)=
\left(\begin{array}{ccc}
0 & -c\,\tau & c_x \\ c\,\tau & 0 & \left(\frac{c_{xx}-c\tau^2}{c}\right) \\ -c_x &  -\left(\frac{c_{xx}-c\tau^2}{c}\right) & 0 
\end{array}\right)
\left(\begin{array}{c}
T\\ n\\ b
\end{array}\right)(t,x).
\end{equation}\\

Now, given a curvature and a torsion obtained by \eqref{Hasimoto} from a solution of \eqref{Hasimoto-eq}, we construct $(T,n,b)$ as explained in subsection \S\ref{constrsection}. We fix an initial condition $(T,n,b)(\tilde{t_0},0)$. Then we define $(T,n,b)(t,0)$ by imposing
\eqref{Frenettimebis} at $x=0$. Finally, $(T,n,b)(t,x)$ is obtained from $(T,n,b)(t,0)$ by integrating the Frenet system for fixed $t$. Showing that $T$  solves indeed 
$$T_t=T\land T_{xx},$$
is then equivalent to showing that 
$$T_t=- c\tau n+c_x b.$$
We shall prove actually that we have the whole system of derivatives in time \eqref{Frenettimebis}. 

Let us introduce the notation $(\alpha,\beta,\gamma) (t,x)$ for those functions such that
$$\left(\begin{array}{c}
T\\ n\\ b
\end{array}\right)_t(t,x)=
\left(\begin{array}{ccc}
0 & \alpha & \beta \\ -\alpha & 0 & \delta \\ -\beta &  -\delta & 0 
\end{array}\right)
\left(\begin{array}{c}
T\\ n\\ b
\end{array}\right)(t,x).$$
By the way we have constructed $T$, it follows that $(\alpha,\beta,\gamma)$ and $(-c\tau, c_x, \frac{c_{xx}-c\tau^2}{c})$ are the same at $x=0$. An easy computation of the derivatives in time and in space of $T$ and of $n$ shows that these functions solve
\begin{equation}\label{last}
\left(\begin{array}{c}
\alpha\\ \beta\\ \gamma
\end{array}\right)_x(t,x)=
\left(\begin{array}{ccc}
0 & \tau & 0 \\ -\tau & 0 & c \\ 0 &  -c & 0 
\end{array}\right)
\left(\begin{array}{c}
\alpha\\ \beta\\ \gamma
\end{array}\right)(t,x)+\left(\begin{array}{c}
c_t\\ 0\\ \tau_t
\end{array}\right)(t,x).
\end{equation}
Let us notice that since $(c,\tau)$ were obtained by \eqref{Hasimoto} from a solution of \eqref{Hasimoto-eq}, they solve DaRios-Betchov's system \cite{DaR}, \cite{Be},
\begin{equation}\label{Betchov}
\left\{\begin{array}{ll}c_t=-2c_x\,\tau-c\,\tau_x,\\ \tau_t=\left(\frac{c_{xx}-c\,\tau^2}{c}\right)_x+ c_x\,c. \end{array}\right.
\end{equation}
A straightforward calculation shows then that $(-c\tau, c_x, \frac{c_{xx}-c\tau^2}{c})$ is also a solution of \eqref{last}. Therefore, for fixed $t$, $(\alpha,\beta,\gamma)$ and $(-c\tau, c_x, \frac{c_{xx}-c\tau^2}{c})$ are two solutions of \eqref{last} with the same initial data at $x=0$. It follows that they coincide for all $(t,x)$, so we obtain the system of derivatives in time \eqref{Frenettimebis}. In particular, we have indeed that $T$ constructed this way solves
$$T_t=T\land T_{xx}.$$

\end{document}